%% file: paper.tex
\documentclass[10pt,a4paper]{amsart} 
\usepackage{graphicx,latexsym,amssymb,color}
\theoremstyle{plain} 
\newtheorem{thm}{Theorem}[section] 
\newtheorem{lem}[thm]{Lemma} 
\newtheorem{prop}[thm]{Proposition} 
\newtheorem{cor}[thm]{Corollary} 
\theoremstyle{definition}
\newtheorem{defi}{Definition} 
\theoremstyle{remark} 
\newtheorem{rem}[thm]{Remark}

\newtheorem*{acknowledgments}{Acknowledgments} 
\numberwithin{equation}{section}
\numberwithin{figure}{section}
\newcommand{\bd}{\begin{description}}   
\newcommand{\ed}{\end{description}} 
\newcommand{\ba}{\begin{array}}      \newcommand{\ea}{\end{array}} 
\newcommand{\bc}{\begin{center}}     \newcommand{\ec}{\end{center}} 
\newcommand{\be}{\begin{enumerate}}  \newcommand{\ee}{\end{enumerate}} 
\newcommand{\beq}{\begin{eqnarray}}  \newcommand{\eeq}{\end{eqnarray}} 
\newcommand{\beQ}{\begin{eqnarray*}} \newcommand{\eeQ}{\end{eqnarray*}} 
\newcommand{\bi}{\begin{itemize}}    \newcommand{\ei}{\end{itemize}}

\newcommand{\ov}{\overline} 
 
\newcommand{\ve}{\varepsilon} 
\newcommand{\s}{\sigma} 
\newcommand{\1}{ {{\mathbf 1}} }
\newcommand{\n}{ \{ 1,...,n \} }
\begin{document} 
\title[]{On $C_{n}$-moves for links} 
\author[J.B. Meilhan]{Jean-Baptiste Meilhan} 
\address{Research Institute for Mathematical Sciences (RIMS)\\
         Kyoto University\\
         Oiwake-cho, Kitashirakawa, Sakyo-ku, \\
         Kyoto 606-8502, Japan}
	 \email{meilhan@kurims.kyoto-u.ac.jp}
\author[A. Yasuhara]{Akira Yasuhara} 
\address{Tokyo Gakugei University\\
         Department of Mathematics\\
         Koganeishi \\
         Tokyo 184-8501, Japan}
	 \email{yasuhara@u-gakugei.ac.jp}

\thanks{The first author is supported by a Postdoctoral
  Fellowship and a Grant-in-Aid for Scientific Research of the Japan
  Society for the Promotion of Science.  
The second author is partially supported by a Grant-in-Aid for Scientific Research (C) 
($\#$18540071) of the Japan Society for the Promotion of Science.}

\date{\today}
\subjclass[2000]{57M25, 57M27}
\keywords{$C_n$-moves, Milnor invariants, string links, Brunnian links, claspers}
\begin{abstract} 
A $C_n$-move is a local move on links defined by Habiro and Goussarov, which can be regarded as a 
`higher order crossing change'.  
We use Milnor invariants with repeating indices to 
provide several classification results for links up to $C_n$-moves, under certain restrictions.  
Namely, we give a classification up to $C_4$-moves of $2$-component links, $3$-component Brunnian 
links and $n$-component $C_3$-trivial links, 
and we classify $n$-component link-homotopically trivial Brunnian links up to $C_{n+1}$-moves.  
\end{abstract} 
\maketitle 
\section{Introduction}
\begin{figure}[!h]
\includegraphics{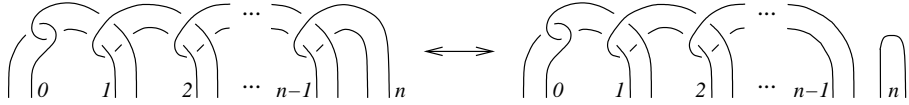}
\caption{A $C_n$-move involves $n+1$ strands of a link, labeled here by integers 
between $0$ and $n$.   } \label{cnm}
\end{figure}

A $C_n$-move is a local move on links as illustrated in Figure \ref{cnm}, 
which can be regarded as a 
kind of `higher order crossing change' (in particular, a $C_1$-move is a crossing change).  
These local moves were introduced by Habiro \cite{Hmt} and independently by Goussarov \cite{G}.   

The $C_n$-move generates an equivalence relation on links, called \emph{$C_n$-equivalence}.  
This notion can also be defined using the theory of claspers (see \S \ref{clasp}). 
The $C_n$-equivalence relation becomes finer as $n$ increases, i.e., the $C_m$-equivalence implies 
the $C_k$-equivalence for $m>k$.  
It is well known that the $C_n$-equivalence allows to approximate the topological 
information carried by Goussarov-Vassiliev invariants.  
Namely, two links cannot be distinguished by any Goussarov-Vassiliev invariant of 
order less than $n$ if
they are $C_n$-equivalent \cite{G,H}.  

Links are classified up to $C_2$-equivalence \cite{MN}.  
$C_3$-classifications of links with $2$ or $3$-components, or of algebraically split links 
are done by Taniyama and the second author \cite{TY}.  
These results involve Milnor $\ov{\mu}$ invariants 
(of length $\le 3$) with distinct indices.  (For the definition of Milnor invariants, 
see \S \ref{milnor}.)  
In this paper, we use Milnor $\ov{\mu}$ invariants with (possibly) repeating indices to 
provide several classification results for higher values of $k$, under certain restrictions.  

First, we consider the case $k=4$.  
We obtain the following for {\em $C_3$-trivial links}, i.e., links which are 
$C_3$-equivalent to the unlink.  

\begin{thm}\label{C4}
Let $L$ and $L'$ be two $n$-component $C_3$-trivial links.  
Then $L$ and $L'$ are $C_4$-equivalent if and only if they satisfy the following properties
 \be
   \item $\ov{\mu}_L(I)=\ov{\mu}_{L'}(I)$ for all multi-index $I$ with $|I|=4$, 
   \item no Vassiliev knot invariant of order $\le 3$ can distinguish 
the $i^{th}$ component of $L$ and the $i^{th}$ 
         component of $L'$, for all $1\le i\le n$.  
 \ee  
\end{thm}

\noindent Here, a multi-index $I$ is a sequence of non-necessarily distinct 
integers in $\{ 1,...,n\}$, 
and $|I|$ denotes the number of entries in $I$.  
Note that this result, together with \cite{MN} and \cite{TY}, imply the 
following.

\begin{cor} 
An $n$-component link $L$ is  
$C_4$-trivial if and only if $\ov{\mu}_L(I)=0$ for all 
multi-index $I$ with $|I|\le 4$ and
any Vassiliev knot invariant of order $\le 3$ vanishes for each component.  
\end{cor}

For $2$-component links, we obtain a refinement of a result of H.A. Miyazawa 
\cite[Thm. 1.5]{M}.  

\begin{prop}\label{C42}
Let $L$ and $L'$ be two $2$-component links.  
Then $L$ and $L'$ are $C_4$-equivalent if and only if they are not 
distinguished by any Vassiliev invariant of order $\le 3$.
\end{prop}

On the other hand, we consider Brunnian links.  
Recall that a link $L$ in the $3$-sphere $S^3$ is \emph{Brunnian} if every proper sublink of $L$ is trivial.  
In particular, all trivial links are Brunnian.  
It is known that an $n$-component link is Brunnian if and only if it can be turned 
into the unlink by a sequence of $C_{n-1}$-moves of a specific type, 
called {\em $C^a_{n-1}$-moves}, involving \emph{all}  
the components \cite{Hb}.  
Also, two $n$-component Brunnian links are $C_n$-equivalent if and only if 
their Milnor invariants $\ov{\mu}(\s(1),...,\s(n-2),n-1,n)$ coincide for 
all $\s$ in the symmetric group $S_{n-2}$ \cite{HM}.  
Here, we consider the next stage, namely $C_{n+1}$-moves for $n$-component Brunnian links.  

Recall that two links are \emph{link-homotopic} if they are related by a sequence of isotopies and 
self-crossing changes, i.e., crossing changes involving two strands of the same component.  
For $n$-component Brunnian links, the link-homotopy coincides with the $C_n$-equivalence \cite{MY,HM}.   
Given $k\in \n$ and a bijection $\tau$ from $\{ 1,...,n-1 \}$ to $\n \setminus \{k \}$, set
  $$ \mu_{\tau}(L) := \ov{\mu}_L(\tau(1),...,\tau(n-1),k,k).  $$ 
We obtain the following.  

\begin{thm}\label{lh}
Two $n$-component link-homotopically trivial Brunnian links $L$ and $L'$ are $C_{n+1}$-equivalent 
if and only if 
$\mu_{\tau}(L)=\mu_{\tau}(L')$ for all $k\in \n$, $\tau\in \mathcal{B}(k)$,   
where $\mathcal{B}(k)$ denotes the set of all bijections $\tau$ from $\{ 1,...,n-1 \}$ to 
$\n \setminus \{k \}$ such that 
$\tau(1)<\tau(n-1)$.  
\end{thm}

In the case of $3$-component Brunnian links, we have the following improvement of Theorem \ref{lh}.  

\begin{thm} \label{trois}
Two $3$-component Brunnian links $L$ and $L'$ are $C_4$-equivalent if and only if 
$\ov{\mu}_L(123)=\ov{\mu}_{L'}(123)$, 
$\ov{\mu}_L(1233)= \ov{\mu}_{L'}(1233)$, 
$\ov{\mu}_L(1322)= \ov{\mu}_{L'}(1322)$ and 
$\ov{\mu}_L(2311)= \ov{\mu}_{L'}(2311)$.\footnote{Note that $\ov{\mu}_L(ijkk)$ denotes here the \emph{residue class} of 
the integer $\mu_L(ijkk)$ (defined in \S \ref{milnor}) modulo $\ov{\mu}_L(ijk)$.  }    
\end{thm}


The rest of the paper is organized as follows.  
In Section \ref{clasp} we recall elementary notions of the theory of claspers.  
In Section \ref{milnor} we recall the definition of Milnor invariants for (string) 
links and give some lemmas.  
In Section \ref{BSL} we consider Brunnian string links.  The main result of 
this section is Proposition \ref{slcn+1}, 
which gives a set of generators for the abelian group of $C_{n+1}$-equivalence 
classes of $n$-component Brunnian string links.  
In Section \ref{BL} we use results of Section \ref{BSL} to prove Theorems \ref{lh} and \ref{trois}.  
In Section \ref{c4eq} we prove Theorem \ref{C4} and Proposition \ref{C42}.  

\begin{acknowledgments}
  The authors wish to thank Kazuo Habiro for helpful comments and conversations.
\end{acknowledgments}
\section{Claspers and local moves on links} \label{clasp}
\subsection{A brief review of clasper theory} \label{review} 
Let us briefly recall from \cite{H} the basic notions of clasper theory for (string) links.  
In this paper, we essentially only need the notion of $C_k$-tree.  
For a general definition of claspers, we refer the reader to \cite{H}.  
\begin{defi}\label{defclasp}
Let $L$ be a link in $S^3$.  
An embedded disk $F$ in $S^3$ is called a {\em tree clasper} for $L$ if 
it satisfies the following (1), (2) and (3):\\
(1) $F$ is decomposed into disks and bands, called {\em edges}, each of which 
connects two distinct disks.\\
(2) The disks have either 1 or 3 incident edges, called {\em leaves} or 
{\em nodes} respectively.\\
(3) $L$ intersects $F$ transversely and the intersections are contained 
in the union of the interior of the leaves. \\
The \emph{degree} of $G$ is the number of the leaves \emph{minus} $1$.  
\end{defi}

A degree $k$ tree clasper is called a $C_k$-tree.  
A $C_k$-tree is \emph{simple} if each leaf intersects $L$ at one point.  

We will make use of the drawing convention for claspers of \cite[Fig. 7]{H}, except 
for the following: a $\oplus$ (resp. $\ominus$) on an edge represents a positive 
(resp. negative) half-twist. (This replaces the 
convention of a circled $S$ (resp. $S^{-1}$) used in \cite{H}).    

\begin{figure}[!h]
\includegraphics{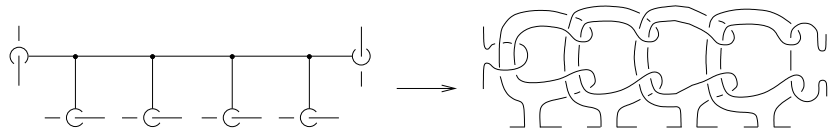}
\caption{Surgery along a simple $C_5$-tree.} \label{ckmove}
\end{figure}

Given a $C_k$-tree $G$ for a link $L$ in $S^3$, there is a procedure to construct, in a regular neighborhood of $G$, a 
framed link $\gamma(G)$. There is thus a notion of \emph{surgery along $G$}, which is defined as surgery 
along $\gamma(G)$.  
There exists a canonical diffeomorphism between $S^3$ and the manifold $S^3_{\gamma(G)}$:   
surgery along the $C_k$-tree $G$ can thus be regarded as a local move on $L$ in $S^3$.  
We say that the resulting link $L_G$ in $S^3$ is obtained by surgery on $L$ along $G$.  
In particular, surgery along a simple $C_k$-tree as illustrated in Figure \ref{ckmove} 
is equivalent to band-summing a copy of the $(k+1)$-component 
Milnor's link $L_{k+1}$ (see \cite[Fig. 7]{Milnor}), 
and is equivalent to a $C_k$-move as defined in the introduction (Figure \ref{cnm}).  
A $C_k$-tree $G$ having the shape of the tree clasper in Figure \ref{ckmove} is called \emph{linear}, 
and the left-most and right-most leaves of $G$ in Figure 2.1 are called the \emph{ends} of $G$.  

The $C_k$-equivalence (as defined in the introduction) coincides with the equivalence relation on links generated 
by surgery along $C_k$-trees and isotopies.  
We use the notation $L\sim_{C_k} L'$ for two $C_k$-equivalent links $L$ and $L'$.  

\subsection{Some lemmas}

In this subsection we give some basic results of calculus of claspers, whose proof can be found in \cite{H} 
or \cite{these}.  
For convenience, we give the statements for string links.  
Recall that a string link is a pure tangle, without closed components (see \cite{HL} for a precise definition).  
Denote by $SL(n)$ the set of $n$-component string links up to isotopy with respect to the 
boundary.  $SL(n)$ has a monoid structure with composition given by the \emph{stacking product}, 
denoted by $\cdot$, and 
with the trivial $n$-component string link $\1_n$ as unit element. 

\begin{lem} \label{cc}
Let $T$ be a union of $C_k$-trees for a string link $L$, and let $T'$ be obtained from $T$ by passing 
an edge across $L$ or 
across another edge of $T$, or by sliding a leaf over a leaf of another component of $T$ (see Figure \ref{lemcc}).  
Then $L_T\sim_{C_{k+1}} L_{T'}$.  
 \begin{figure}[!h]
  \includegraphics{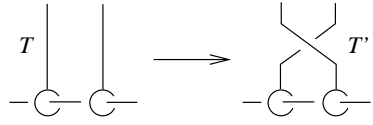}
  \caption{Sliding a leaf over another leaf.  }\label{lemcc}
 \end{figure}
\end{lem}

\begin{lem} \label{twist}
 Let $T$ be a $C_k$-tree for $\1_n$ and let $\ov{T}$ be a $C_k$-trees 
obtained from $T$ by adding a half-twist on an edge. Then 
 $(\1_n)_T\cdot (\1_n)_{\ov{T}}\sim_{C_{k+1}} \1_n$.  
\end{lem}

\begin{lem} \label{asihx}
Consider some $C_k$-trees $T$ and $T'$ (resp. $T_I$, $T_H$ and $T_X$) for $\1_n$ which differ only in a small ball 
as depicted in Figure \ref{F05}, then $(\1_n)_{T}\cdot (\1_n)_{T'}\sim_{C_{k+1}} \1_n$  
(resp. $(\1_n)_{T_I}\sim_{C_{k+1}} (\1_n)_{T_H}\cdot (\1_n)_{T_X}$). 
 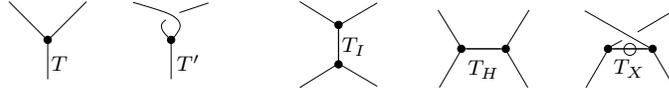
\begin{figure}[!h]
  \input{asihx.pstex_t}
  \caption{The AS and IHX relations for $C_k$-trees.  }\label{F05}
 \end{figure}
\end{lem} 

\begin{lem} \label{split}
Let $G$ be a $C_k$-tree for $\1_n$. Let $f_1$ and $f_2$ be two disks obtained by splitting a leaf $f$ of $G$ 
along an arc $\alpha$ as shown in figure \ref{scind} 
(i.e., $f=f_1\cup f_2$ and $f_1\cap f_2=\alpha$). 
Then, $(\1_n)_G\sim_{C_{k+1}} (\1_n)_{G_1}\cdot (\1_n)_{G_2}$, 
where $G_i$ denotes the $C_k$-tree for $\1_n$ obtained from $G$ by replacing $f$ by $f_i$ ($i=1,2$).  
\end{lem}
 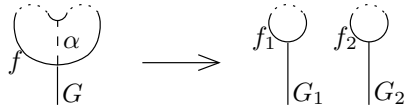
\begin{figure}[!h]
  \input{split.pstex_t}
  \caption{Splitting a leaf.  }\label{scind}
 \end{figure}

\subsection{$C^a_k$-trees and $C^a_k$-equivalence}\label{sec:cak-equivalence}

\begin{defi} \label{r9}
  Let $L$ be an $m$-component link in a $3$-manifold $M$.  For $k\ge m-1$,
  a (simple) $C_k$-tree $T$ for $L$ in $M$ is a \emph{(simple) $C^a_k$-tree} if it satisfies the following: 
  \begin{enumerate}
  \item For each disk-leaf $f$ of $T$, $f\cap L$ is contained in a single component of $L$, 
  \item $T$ intersects \emph{all} the components of $L$.
  \end{enumerate}
\end{defi}

The \emph{$C^a_k$-equivalence} is an equivalence relation on links generated by surgeries along $C^a_k$-trees 
and isotopies.  The next result shows the relevance of this notion in the study of Brunnian (string) links. 

\begin{thm}[\cite{Hb,MY}] \label{brunian}
 Let $L$ be an $n$-component link in $S^3$.  Then $L$ is Brunnian if and only if 
 it is $C^a_{n-1}$-equivalent to the $n$-component trivial link.  
\end{thm}

Further, it is known that, for $n$-component Brunnian links the $C_n$-equivalence coincides with the 
$C^a_n$-equivalence (and with the link-homotopy) \cite{MY}, see also \cite{HM}.  

\section{On Milnor invariants} \label{milnor}

\subsection{A short definition} 

J. Milnor defined in \cite{Milnor} a family of invariants of oriented, ordered 
links in $S^3$, known as Milnor's $\ov{\mu}$-invariants. 

Given an $n$-component link $L$ in $S^3$, denote by $\pi$ the fundamental group of $S^3\setminus L$, and by $\pi_q$ 
the $q^{th}$ subgroup of the lower central series of $\pi$.  
We have a presentation  of $\pi/ \pi_q$ with $n$ generators, given by a meridian $m_i$ of the $i^{th}$ component of $L$.  
So for $1\le i\le n$, the longitude $l_i$ of the $i^{th}$ component of $L$ is expressed modulo $\pi_q$ as a word  
in the $m_i$'s (abusing notations, we still denote this word by $l_i$).  

The \emph{Magnus expansion} $E(l_i)$ of $l_i$ is the formal power series in 
non-commuting variables $X_1,...,X_n$ obtained by 
substituting $1+X_j$ for $m_j$ and $1-X_j+X_j^2-X_j^3+...$ for $m_j^{-1}$, $1\le j\le n$.  
We use the notation $E_k(l_i)$ to denote the degree $k$ part of $E(l_i)$ 
(where the degree of a monomial in the $X_j$ is simply 
defined by the sum of the powers).  

Let $I=i_1 i_2 ...i_{k-1} j$ be a multi-index (i.e., a sequence of possibly repeating 
indices) among $\n$. 
Denote by $\mu_L(I)$ the coefficient of $X_{i_1}...X_{i_{k-1}}$ in the Magnus expansion $E(l_j)$.  
\emph{Milnor invariant} $\ov{\mu}_L(I)$ is the residue class of $\mu_L(I)$ modulo the greatest common divisor of 
all Milnor invariants $\mu_L(J)$ such that $J$ is obtained from $I$ by removing at least one index and permuting 
the remaining indices cyclicly.  
$|I|=k$ is called the \emph{length} of Milnor invariant $\ov{\mu}_L(I)$.  
 
The indeterminacy comes from the choice of the meridians $m_i$.  Equivalently, it comes from the indeterminacy of 
representing the link as the closure of a string link \cite{HL}.  
Indeed, $\mu(I)$ is a well-defined invariant for string links.  
Furthermore, $\mu(I)$ is known to be a Goussarov-Vassiliev invariant of degree $|I|-1$ for string 
links \cite{BN2,Lin}.  


\subsection{Some lemmas}

Let us first recall a result due to Habiro.

\begin{lem}[\cite{H}] \label{milnCk}
 Milnor invariants of length $k$ for (string) links are invariants of $C_{k}$-equivalence.  
\end{lem}

Next we state a simple lemma that will be used in the following.
\begin{lem} \label{miln0}
 Let $L$ be an $n$-component string link which is obtained from $\1_n$ by surgery along 
a union $F$ of $C_k$-trees which 
 is disjoint from the $j^{th}$ component of $\1_n$.  
 Then $\mu_L(I)=0$, for all multi-index $I$ containing $j$ and satisfying $|I|\le k+1$.  
\end{lem}

\begin{proof}
Consider a diagram of $\1_n$ together with $F$.  The diagram contains several crossings between 
an edge of $F$ and the $j^{th}$ component of $\1_n$.  Denote by $F_o$ (resp. $F_u$) 
the union of $C_k$-trees 
obtained from $F$ by performing crossing changes so that the $j^{th}$ component of 
$\1_n$ overpasses (resp. underpasses) all edges.  
By Lemma \ref{cc}, we have $L\sim_{C_{k+1}} U_{F_o}\sim_{C_{k+1}} U_{F_u}$.  
The result then follows from Lemma \ref{milnCk} and the following observation.  

Consider the diagram $D$ of a string link $K$.  If the $i^{th}$ component of 
$K$ overpasses all the other components 
in $D$, it follows from the definition of Milnor invariants that $\mu_K(I)=0$ 
for any sequence $I$ with the last index $i$.  
Similarly, if the $i^{th}$ component of $K$ underpasses all the other components in 
$D$, then $\mu_K(I)=0$ for any sequence 
$I$ with an index $i$ which is not equal to the last one. 
\end{proof}

We have the following simple additivity property.  

\begin{lem}\label{add}
Let $L$ and $L'$ be two $n$-component string links such that all Milnor invariants of $L$ 
(resp. $L'$) of length $\le m$ 
(resp. $\le m'$) vanish.  
Then $\mu_{L\cdot L'}(I)=\mu_{L}(I)+\mu_{L'}(I)$ for all $I$ of length $\leq m+m'$.  
\end{lem}

\begin{proof}
Milnor invariant of $L\cdot L'$ is computed by taking the Magnus expansion of the 
$k^{th}$ longitude $L_k$ of $L\cdot L'$.  
Denote respectively by $l_i$ and $m_i$ (resp. $l'_i$ and $m'_i$)  
the $i^{th}$ meridian and longitude of $L$ (resp. $L'$) ; $1\le i\le n$.  
We have $L_k=l_k\cdot \tilde{l}'_k$, 
where $\tilde{l}'_k$ is obtained from $l'_k$ by replacing $m'_i$ with 
$M_i=l_i^{-1} m_i l_i$ for each $i$.  
So $E(L_k)=E(l_k)\cdot E(\tilde{l}'_k)$, where 
$E(\tilde{l}'_k)$ is obtained from $E(l'_k)$ by substituting $\tilde{X}_i$ for $X_i$ in $E(l'_k)$, 
where $\tilde{X_i}:=E(M_i)-1$.  

The Magnus expansion of $l_i$ is the form $E(l_i)= 1 +($terms of order $\ge m)$, so 
\beQ
E(M_i) & = & E(l_i^{-1}) E(m_i) E(l_i) \\
 & = & E(l_i^{-1})E(l_i) + E(l_i^{-1})X_i E(l_i) \\
 & = & 1 + X_i + \textrm{ (terms of order $\ge m+1$)}.   
\eeQ

So 
$E(\tilde{l}'_k)$ is obtained from $E(l'_k)=\sum_{j\ge m'} E_{j}(l'_k)$ 
by replacing each $X_i$ by $X_i +$ (terms of order $\ge m+1$) for all $i$.  It follows that  
  $$ E(\tilde{l}'_k) = 1 + \sum_{m+m'-1\ge j\ge m'} E_{j}(l'_k) + \left(\textrm{ terms of 
order $\ge (m+m')$}\right). $$

It follows that $ E(L_k) = E(l_k) E(\tilde{l}'_k)$ has the form  
  $$ 1 + \sum_{m+m'-1\ge j\ge m}E_{j}(l_k) + \sum_{m+m'-1\ge j\ge m'}E_{j}(l'_k) +
\left(\textrm{ terms of order $\ge (m+m')$}\right), $$
which implies that all Milnor invariants of length $\le m+m'$ of $L\cdot L'$ are additive.  
\end{proof}

\section{$C_{n+1}$-moves for $n$-component Brunnian string links} \label{BSL}

An $n$-component string link $L$ is Brunnian if every proper substring link of $L$ is the trivial string link.  
In particular, any trivial string link is Brunnian.  
$n$-component Brunnian string links form a submonoid of $SL(n)$, denoted by $BSL(n)$.  

Recall that, given $L\in SL(n)$, the \emph{closure} $\mathrm{cl}(L)$ of $L$ is an $n$-component 
link in $S^3$ \cite{HL}.  
By \cite{Hb}, an $n$-component link is Brunnian if and only if it is the closure of a 
certain Brunnian string link.  

\subsection{$n$-component Brunnian string links up to $C_n$-equivalence}  \label{slcn}

Let ${BSL}(n)/C_n$ denote the abelian group of $C_n$-equivalence classes of $n$-component 
Brunnian string links.  
Habiro and the first author gave in \cite{HM} a basis for ${BSL}(n)/C_n$ as follows.  

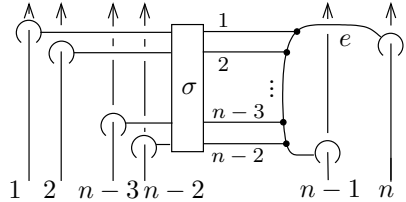
\begin{figure}[!h]
  \input{Tsigm.pstex_t}
  \caption{The simple $C^a_n$-tree $T_\s$.  Here, the numbering of the edges just indicates how 
  $\s\in S_{n-1}$ acts on the edges of $T_\s$ (a similar notation is used in Fig. \ref{Ttau}).  }\label{Tsigma}
 \end{figure}

Let $\s$ be an element in the symmetric group $S_{n-2}$.  Denote by $L_\s$ the $n$-component string link obtained 
from $\1_n$ by surgery along the $C^a_{n-1}$-tree $T_\s$ shown in Figure \ref{Tsigma}.  
Likewise, denote by $(L_\s)^{-1}$ the $n$-component string link obtained from the $C^a_{n-1}$-tree $\ov{T}_\s$,  
which is obtained from $T_\s$ by adding a positive half-twist in the edge $e$  (see Figure \ref{Tsigma}).  
 
Let $\mu_\s(L)$ denote the Milnor invariant $\mu_L(\s(1),...,\s(n-2),n-1,n)$ for any element $\s\in S_{n-2}$.  

\begin{prop}[\cite{HM}] \label{hm}
Let $L$ be an $n$-component Brunnian string link.  Then 
 $$ L\sim_{C_n} \prod_{\s \in S_{n-2}} (L_{\s})^{\mu_{\s}(L)}. $$  
\end{prop}

\begin{rem} \label{remHM}
Recall from \cite{HM,MY} that the $C_n$-equivalence, the link-homotopy and the $C^a_n$-equivalence 
coincide on $BSL(n)$.   
\end{rem}

\subsection{$n$-component Brunnian string links up to $C_{n+1}$-equivalence} 

In this section, we study the quotient $BSL(n) / C_{n+1}$.  Note that $BSL(n) / C_{n+1}$ 
is a finitely generated abelian group 
(this is shown by using the same arguments as in the proof of \cite[Lem. 5.5]{H}).  

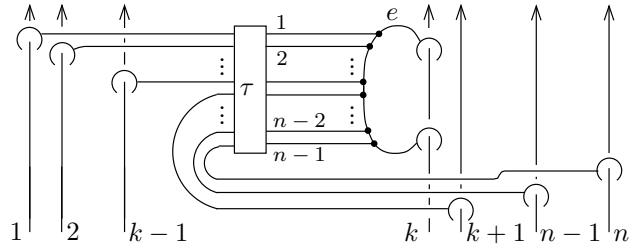
\begin{figure}[!h]
  \input{Tta.pstex_t}
  \caption{The simple $C^a_n$-tree $G_\tau$.  }\label{Ttau}
 \end{figure}

Given $k\in \n$, consider a bijection $\tau$ from $\{ 1,...,n-1 \}$ to $\n \setminus \{k \}$.  
Denote by $V_\tau$ the $n$-component string link obtained from $\1_n$ by surgery along the 
$C^a_{n}$-tree $G_\tau$ shown in Figure \ref{Ttau}.  
Denote by $\ov{G}_\tau$ the $C^a_{n}$-tree for $\1_n$ obtained from $G_\tau$ by adding a 
positive half-twist in the edge $e$  (see Figure \ref{Tsigma}).  
Let $(V_\tau)^{-1}$ be the $n$-component string link obtained from $\1_n$ by 
surgery along $\ov{G}_\tau$.    

Set
  \begin{equation*}
    \mu_{\tau}(L) := \mu_L(\tau(1),...,\tau(n-1),k,k) 
  \end{equation*}  
Denote by $\mathcal{B}(k)$ the set of all bijections $\tau$ from $\{ 1,...,n-1 \}$ to 
$\n \setminus \{k \}$ such that $\tau(1)<\tau(n-1)$, and 
denote by $\rho$ a bijection from $\{ 1,...,n-1 \}$ to itself defined by $\rho(i)=n-i$.  
We have the following. 

\begin{prop}\label{slcn+1}
Let $L$ be an $n$-component Brunnian string link.  Then 
\begin{equation} \label{1}
  L\sim_{C_{n+1}} \left(\prod_{\s \in S_{n-2}} (L_{\s})^{\mu_{\s}(L)}\right)
\cdot L_1\cdot ... \cdot L_n,
\end{equation}
where, for each $k~(1\le k\le n)$, $L_k$ is the $n$-component Brunnian string link
  $$ \prod_{\tau\in \mathcal{B}(k)} (V_{\tau})^{n_{\tau}(L)}\cdot (V_{\tau \rho})^{n'_{\tau}(L)}, $$
such that, for any $\tau\in \mathcal{B}(k)$ ($k=1,...,n$), $n_{\tau}(L)$ and $n'_{\tau}(L)$ are two integers satisfying 
\begin{equation} \label{ega}
  n_{\tau}(L)+(-1)^{n-1}n'_{\tau}(L)=\mu_{\tau}(L_1\cdot ... \cdot L_n). 
\end{equation}
\end{prop}

\begin{proof}
By Proposition \ref{hm} and Remark \ref{remHM}, $L$ is obtained from the $n$-component string link 
 $$ L_0:=\prod_{\s \in S_{n-2}} (L_{\s})^{\mu_{\s}(L)} $$ 
by surgery along a disjoint union $F$ of simple $C^a_n$-trees.  
By Lemma \ref{cc}, we have
 $$ L \sim_{C_{n+1}} L_0\cdot (\1_n)_{G_1}\cdot ...\cdot (\1_n)_{G_p}, $$
where, $G_j$ $(1\le j\le p)$ are simple $C^a_n$-trees for $\1_n$.  
Denote by $k_j$ the (unique) element of $\n$ such that $G_j$ intersects twice the $k_j^{th}$ component 
of $\1_n$ ($1\le j\le p$). 
We can use the AS and IHX relations for tree claspers to replace, up to $C_{n+1}$-equivalence, 
each of these $C^a_n$-trees with a union of linear $C^a_n$-trees whose ends intersect 
the $k_j^{th}$ component.  
More precisely, by lemmas \ref{asihx}, \ref{twist} and \ref{cc} we have for each $1\le j\le p$
 $$ (\1_n)_{G_j} \sim_{C_{n+1}} \prod_{i=1}^{m_j} (V_{\nu_{ij}})^{\ve_{ij}}, $$	
where $\ve_{ij}\in \mathbf{Z}$ and where $\nu_{ij}$ is a bijection from 
$\{ 1,...,n-1 \}$ to $\n\setminus \{ k_j\}$.  
Since there exists, for each such $\nu_{ij}$, a unique element $\tau$ of $\mathcal{B}(k_j)$ 
such that $\nu_{ij}$ is  equal 
to either $\tau$ or $\tau \rho$, it follows 
that $L$ is $C_{n+1}$-equivalent to an $n$-component string link of the form given in (\ref{1}).  
It remains to prove (\ref{ega}).  

First, let us compute $\mu_{\tau}(V_{\eta})$ for all $\tau\in \mathcal{B}(k)$ 
and $\eta\in \mathcal{B}(l)$ ; $k,l=1,...,n$.  
By \cite[Theorem 7]{Milnor2}, we have 
  $$ \mu_{\tau}(V_{\eta}) = \mu_{\tau,n+1}(W_{\eta}), $$
where $\mu_{\tau,n+1}$ is Milnor invariant 
$\mu(\tau(1),...,\tau(k-1),\tau(k+1),...,\tau(n),k,n+1)$ and where $W_{\eta}$ 
denotes the $(n+1)$-component string link obtained from $V_{\eta}$ by taking, as the $(n+1)^{th}$ 
component, a parallel copy of the $k^{th}$ 
component (so that the $k^{th}$ and the $(n+1)^{th}$ components of $W_{\eta}$ have linking number zero).  
Now recall that $V_{\eta}\cong (\1_n)_{G_{\eta}}$, where $G_{\eta}$ is a $C^a_{n}$-tree as 
shown in Figure \ref{Ttau}.  
So $W_{\eta}\cong (\1_{n+1})_{\tilde{G}_{\eta}}$, where $\tilde{G}_{\eta}$ is a 
$C^a_n$-tree obtained from $G_{\eta}$ by 
replacing each leaf intersecting the $k^{th}$ component of $\1_n$ with a leaf intersecting 
components $k$ and $n+1$ as depicted 
in Figures \ref{F02} and \ref{F01}.  

If $k\ne l$, then $\tilde{G}_{\eta}$ contains exactly one leaf $f$ intersecting 
both the $k^{th}$ and the $(n+1)^{th}$ 
components of $\1_{n+1}$.  
\begin{figure}[!h]
  \input{F02.pstex_t}
  \caption{Here (and in subsequent figures) we fix, for simplicity, 
 	   $n=4$, $k=1$, $l=4$ and $\eta$ is the cyclic permutation $(231)\in S_3$}\label{F02}
 \end{figure}
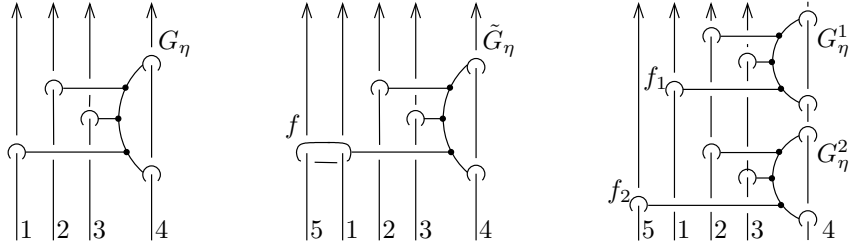
 
\noindent By Lemma \ref{split}, we have 
  $$ (\1_{n+1})_{\tilde{G}_{\eta}}\sim_{C_{n+1}} (\1_{n+1})_{G^1_{\eta}}\cdot (\1_{n+1})_{G^2_{\eta}}, $$
where $G^i_{\eta}$ denotes the simple $C_n$-tree for $\1_{n+1}$ 
obtained from $\tilde{G}_{\eta}$ by replacing $f$ by 
$f_i$ as shown in Figure \ref{F02} ($i=1,2$).  
By Lemmas \ref{milnCk} and \ref{add}, $\mu_{\tau}(V_{\eta})$ is thus equal to 
$\mu_{\tau,n+1}((\1_{n+1})_{G^1_{\eta}}) + \mu_{\tau,n+1}((\1_{n+1})_{G^2_{\eta}})$.  
It follows from Lemma \ref{miln0} that $\mu_{\tau}(V_{\eta})=0$.  
 
Now suppose that $k=l$. Then $\tilde{G}_{\eta}$ contains two leaves intersecting both 
the $k^{th}$ and the $(n+1)^{th}$ components of 
$\1_{n+1}$.  By Lemma \ref{split}, we obtain  
 $$ (\1_{n+1})_{\tilde{G}_{\eta}}\sim_{C_{n+1}} (\1_{n+1})_{G^1_{\eta}}\cdot 
    (\1_{n+1})_{G^2_{\eta}}\cdot (\1_{n+1})_{G^3_{\eta}}\cdot (\1_{n+1})_{G^4_{\eta}}, $$
where, for $1\le i\le 4$, $G^i_{\eta}$ is a simple $C_n$-tree for $\1_{n+1}$ as depicted in Figure \ref{F01}.  
 \begin{figure}[!h]
  \input{F01.pstex_t}
  \caption{}\label{F01}
 \end{figure}
 
\noindent By Lemmas \ref{milnCk}, \ref{miln0} and \ref{add}, it follows that 
 $$\mu_{\tau}(V_{\eta})=\mu_{\tau,n+1}((\1_{n+1})_{G^3_{\eta}})+\mu_{\tau,n+1}((\1_{n+1})_{G^4_{\eta}}). $$  
Observe that the closure of each of these two string links is a copy of Milnor's link \cite[Fig. 7]{Milnor}.  
By a formula of Milnor \cite[pp. 190]{Milnor}, we obtain 
$\mu_{\tau,n+1}((\1_{n+1})_{G^3_{\eta}})=\delta_{\tau,\eta}$, and 
$\mu_{\tau,n+1}((\1_{n+1})_{G^4_{\eta}})=0$, 
where $\delta$ denotes Kronecker's symbol.  
So we obtain that 
 $$ \mu_{\tau}(V_{\eta})=\delta_{\tau,\eta}. $$
Moreover, it follows from Lemmas \ref{add} and \ref{twist} that 
$\mu_{\tau}((V_{\eta})^{-1})=-\delta_{\tau,\eta}$.  

Now consider the string link $V_{\eta \rho}$.  By the same arguments as above, we have that 
$\mu_{\tau}(V_{\eta \rho})=\mu_{\tau}((V_{\eta \rho})^{-1})=0$ if $k\ne l$.  
If $k=l$, it follows from the same arguments as above that 
 $$ \mu_{\tau}(V_{\eta \rho})=\mu_{\tau,n+1}((\1_{n+1})_{G^1_{\eta \rho}})+
\mu_{\tau,n+1}((\1_{n+1})_{G^2_{\eta \rho}}), $$
where $G^1_{\eta \rho}$ and $G^2_{\eta \rho}$ are two simple $C^a_n$-trees for $\1_{n+1}$ as 
depicted in Figure \ref{F03}.  
 \begin{figure}[!h]
  \input{F03.pstex_t}
  \caption{}\label{F03}
 \end{figure}
 
\noindent By Lemma \ref{cc} and isotopy, $(\1_{n+1})_{G^i_{\eta \rho}}$ is 
$C_{k+1}$-equivalent to $(\1_{n+1})_{T^i_{\eta}}$, 
where $T^i_{\eta}$ is as shown in Figure \ref{F03}, $i=1,2$.  
By Lemma \ref{twist}, we thus obtain 
 $$ \mu_{\tau}(V_{\eta \rho})=(-1)^{n-1} \delta_{\tau,\eta}. $$ 

We conclude that 
 $$ \mu_{\tau}(L_1\cdot ... \cdot L_p)=\sum_{1\le i\le p} \mu_{\tau}(L_i)=n_{\tau}(L)+(-1)^{n-1}n'_{\tau}(L). $$ 
\end{proof}

\begin{rem} \label{remlink}
Observe that we obtain the following as a byproduct of the proof of Proposition \ref{slcn+1}.  
Consider the $n$-component Brunnian link $B_{\tau}$ represented in Figure \ref{Btau}, 
for some $\tau\in \mathcal{B}(k)$.  
$B_{\tau}$ is the closure of the $n$-component string link $V_{\tau}$ considered above. 
We showed that, for $1\le l\le n$ and $\eta \in \mathcal{B}(l)$, 
  \begin{equation*}
    \ov{\mu}_{\eta} (B_{\tau})={\mu}_{\eta} (B_{\tau})= \delta_{\eta,\tau}.  
  \end{equation*}
 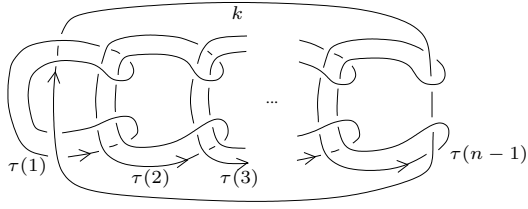
\begin{figure}[!h]
  \input{link.pstex_t}
  \caption{The link $B_{\tau}$.  }\label{Btau}
 \end{figure}
\end{rem}

We conclude this section by showing that 
the string links $V_\tau$ and $V_{\tau \rho}$ are linearly independent in $BSL(n) / C_{n+1}$.  

\begin{prop} 
For any integer $k$ in $\{1,...,n\}~(n\ge 3)$ and any element $\tau\in \mathcal{B}(k)$, 
we have $V_\tau\nsim_{C_{n+1}} V_{\tau \rho}$ nor $(V_{\tau \rho})^{-1}$.
\end{prop}

\begin{rem}
In contrast to the lemma above, we will see in the proof of Proposition \ref{clos} 
that
$\mathrm{cl}(V_\tau)\sim_{C_{n+1}} \mathrm{cl}(V_{\tau \rho})$ or 
$\mathrm{cl}((V_{\tau \rho})^{-1})$.
\end{rem}

\begin{proof}
Consider a diagram of an $n$-component string link $L$.  $L$ lives in a copy of $D^2\times I$ 
standardly embedded in $S^3$.  
The \emph{origin} (resp. \emph{terminal}) of the $i^{th}$ component of $L$ is 
the starting point (resp. ending point) 
of the component, according to the orientation of $L$.  We can construct a knot $K_\tau(L)$ in $S^3$ as follows.  

Connect the terminals of the $k^{th}$ and the ${\tau(1)}^{th}$ components by an arc $a_1$ 
in $S^3\setminus (D^2\times I)$.  
Next, connect the origins of the ${\tau(1)}^{th}$ and the ${\tau(2)}^{th}$ components 
by an arc $a_2$ in $S^3\setminus (D^2\times I)$ 
disjoint from $a_1$, then the terminals of the ${\tau(2)}^{th}$ and the ${\tau(3)}^{th}$ components 
by an arc $a_3$ in $S^3\setminus (D^2\times I)$ disjoint from 
$a_1\cup a_2$.  Repeat this construction until reaching the last component, 
the ${\tau(n-1)}^{th}$ component, and connect the terminal or the origin 
(depending on whether $n$ is even or odd) to the origin of the $k^{th}$ component 
by an arc $a_n$ in $S^3\setminus (D^2\times I)$ 
disjoint from $\bigcup_{1\le i\le n-1} a_i$.  The arcs are chosen so that, if $a_i$ and $a_j$ ($i<j$) meet in the 
diagram of $L$, then $a_i$ overpasses $a_j$.  The orientation of $K_\tau$ is the one induced 
from the $k^{th}$ component.  
An example is given in Figure \ref{Ktau} for the case $n=4$, $k=4$ and $\tau=(231)\in S_3$.   
 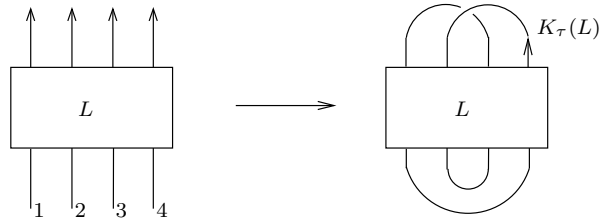
\begin{figure}[!h]
  \input{Ktau.pstex_t}
  \caption{The knot $K_{\tau}(L)$.  }\label{Ktau}
 \end{figure}

It follows immediately from the above construction and \cite[Thm. 1.4]{horiuchi} that 
 $$ P^{(n)}_0(K_\tau(V_\tau);1)=\pm n! 2^n\textrm{ and } P^{(n)}_0(K_\tau(V_{\tau\rho});1)=
P^{(n)}_0(K_\tau((V_{\tau\rho})^{-1});1)=0, $$
where $P^{(k)}_l(K;1)$ denotes the $k^{th}$ derivative of the coefficient 
polynomial $P_k(K;t)$ of $z^k$ in the HOMFLY polynomial 
$P(K;t,z)$ of a link $K$, evaluated in $1$.  
The result then follows from \cite[Cor. 6.8]{H} and the fact that $P^{(n)}_0(K;1)$ is a 
Goussarov-Vassiliev invariant of degree 
$\le n$ \cite{KM}.  
\end{proof}

\section{$C_{n+1}$-moves for $n$-component Brunnian links} \label{BL}

In this section we prove Theorems \ref{lh} and \ref{trois}.  
Let us begin with stating the following link version of Proposition \ref{slcn+1}.  

\begin{prop} \label{clos}
Let $L$ be an $n$-component Brunnian link.  Then 
  $$ L\sim_{C_{n+1}} \mathrm{cl}\left( \prod_{\s \in S_{n-2}} (T_{\s})^{\mu_{\s}(L)}\cdot 
\prod_{1\le k\le n} L'_k \right), $$ 
where, for each $i~(1\le i\le n)$, 
  $$ L'_k:=\prod_{\tau_k\in \mathcal{B}(k)} (V_{\tau_k})^{\mu_{\tau}(L'_1\cdot ...\cdot L'_p)}. $$
\end{prop}

\begin{proof}
By Proposition \ref{slcn+1}, $L$ is $C_{n+1}$-equivalent to the closure of the string link 
\begin{equation}  \label{Ls}
   l=\prod_{\s \in S_{n-2}} ((\1_n)_{T_\s})^{\mu_{\s}(L)}\cdot \prod_{1\le k\le n} 
    \prod_{\tau\in \mathcal{B}(k)} ((\1_n)_{G_{\tau}})^{n_{\tau}(L)}\cdot ((\1_n)_{G_{\tau \rho}})^{n'_{\tau}(L)},   
\end{equation}
where $n_{\tau}(L)$ and $n'_{\tau}(L)$ are two integers satisfying (\ref{ega}).  
Denote by $F$ the union of all the tree claspers involved in (\ref{Ls}), that is $l=(\1_n)_F$.  

For some $k\in \n$, $\tau\in \mathcal{B}(k)$, let $G$ be a copy of the simple $C_n$-tree $G_{\tau \rho}$ in $F$.  
Let $f$ be a leaf of $G$ which intersects $k^{th}$ component of $\1_n$ (Figure \ref{F04}).  
When we close the $k^{th}$ component of $\1_n$, we can slide $f$ over leaves of 
the components of $F\setminus G$ until we 
obtain the $C_n$-tree $G'$ of Figure \ref{F04}.  Denote by $F'$ the union of tree claspers 
obtained from $F$ by this operation.  
By Lemma \ref{cc}, we have $\mathrm{cl}((\1_n)_F)\sim_{C_{n+1}} \mathrm{cl}((\1_n)_{F'})$.  
 \begin{figure}[!h]
  \input{F04.pstex_t}
  \caption{}\label{F04}
 \end{figure}
 
\noindent By Lemma \ref{cc} and isotopy, $(\1_{n})_{G'}$ is $C_{n+1}$-equivalent to $(\1_{n})_{G''}$, 
where $G''$ is the $C_n$-tree depicted in Figure \ref{F04}.  
$G''$ differs from a copy of $G_{\tau}$ by $(n+1)$ half-twists on its edges.  
It thus follows from Lemma \ref{twist} that 
  \begin{equation*}
     \mathrm{cl}((\1_n)_{G_{\tau}}\cdot (\1_n)_{G_{\tau \rho}})\sim_{C_{n+1}}  \left\{ \begin{array}{ll}
      \mathrm{cl}(\1_n) & \text{if $n$ is even,} \\
      \mathrm{cl}(((\1_n)_{G_{\tau}})^2) & \text{if $n$ is odd.}
    \end{array}\right.
  \end{equation*}
$L$ is thus $C_{n+1}$-equivalent to the closure of the string link 
\begin{equation*}  
   \prod_{\s \in S_{n-2}} ((\1_n)_{T_{\s}})^{\mu_{\s}(L)}\cdot \prod_{1\le k\le n} 
    \prod_{\tau\in \mathcal{B}(k)} ((\1_n)_{G_{\tau}})^{n_{\tau}(L)+(-1)^{n-1}n'_{\tau}(L)}.   
\end{equation*}
The result follows from (\ref{ega}).  
\end{proof}

\subsection{The link-homotopically trivial links case: Proof of Theorem \ref{lh}}

\begin{proof}[Proof of Theorem \ref{lh}]
By Proposition \ref{hm}, if an $n$-component Brunnian link $B$ is link-homotopically trivial, 
then $\mu_\s(B)=0$ for all $\s\in S_{n-2}$.  
For all $\tau\in \mathcal{B}(k)$, $k=1,...,n$, ${\mu}_\tau(B)$ is thus a well-defined integer, 
which satisfies ${\mu}_\tau(B)=\mu_\tau(L(B))$ for any string link $L(B)$ whose closure is $B$.  
By Proposition \ref{clos}, we have 
 $$ B\sim_{C_{n+1}} \mathrm{cl}\left( \prod_{1\le k\le n} 
\prod_{\tau\in \mathcal{B}(k)} (V_{\tau})^{\mu_{\tau}(B)}\right) . $$
The result follows immediately.  
\end{proof}

\subsection{The $3$-component links case: Proof of Theorem \ref{trois}} \label{S3}

\begin{proof}[Proof of Theorem \ref{trois}]
The `only if' part of Theorem \ref{trois} is an immediate consequence of Lemma \ref{milnCk}.  
Here we prove the `if' part.  

Let $L$ be a $3$-component Brunnian link.  
By Proposition \ref{clos}, we have 
\begin{equation}\label{e1}
 L\sim_{C_4} \mathrm{cl}(L_0\cdot L_1\cdot L_2\cdot L_3),\quad \textrm{with}\quad L_0:=B^{{\mu}_L(123)}
\textrm{ and } L_p:=V_p^{n_p}~(p=1,2,3),
\end{equation}
where $B$ and $V_p$ ($p=1,2,3$) are $3$-component string links obtained from $\1_3$ by 
surgery along a $C_2$-tree and $C_3$-trees as shown 
in Figure \ref{bor} respectively, and where $n_k=\mu_{L_1\cdot L_2\cdot L_3}(ijkk)$ with $\{i,j,k\}=\{1,2,3\}$ and $i<j$.  Note that ${\mu}_L(123)=\ov{\mu}_L(123)$ since $L$ is Brunnian. 
 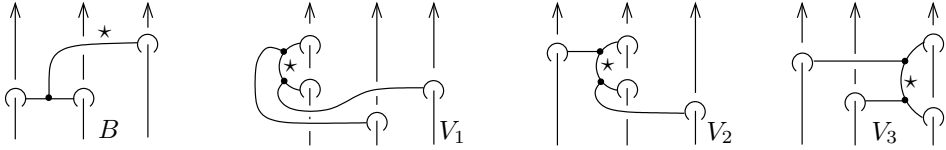
\begin{figure}[!h]
  \input{bor.pstex_t}
  \caption{Here $B^{-1}$ (resp. $V_{p}^{-1}$, $1\le p\le 3$) is defined as obtained from $B$ (resp. 
           $V_{p}$, $1\le p\le 3$) by a positive half-twist on the edge marked by a $\star$.  }\label{bor}
 \end{figure}

We now make an observation.  
Consider a union $Y$ of $k$ parallel copies of a simple $C_2^a$-tree for the $3$-component 
unlink $U=U_1\cup U_2\cup U_3$, 
and perform an isotopy as illustrated in Figure \ref{iso}.  Denote by $Y'$ the resulting union of $C_2$-trees.  
By \cite[Prop. 4.5]{H}, $Y'$ can be deformed into $Y$ by a sequence of $k$ $C_3$-moves, corresponding to $k$ parallel 
copies of a simple $C_3$-tree intersecting twice $U_i$ and once $U_j$ and $U_k$.  
 \begin{figure}[!h]
  \input{isotopy.pstex_t}
  \caption{}\label{iso}
 \end{figure}
So by Lemma \ref{twist}, $U_Y$ is $C_4$-equivalent to 
$\mathrm{cl}\left( (\1_n)_Y\cdot (\1_n)^{\pm k}_{V_i} \right)$.\footnote{Here, 
abusing notations, we still denote by $Y$ a union of $k$ simple $C_2$-trees for $\1_3$ 
such that $\mathrm{cl}((\1_3)_Y)\cong U_Y$.  } 
Note that for any union $F$ of $C_3$-trees, 
$U_{Y\cup F}\sim_{C_4}\mathrm{cl}((\1_n)_{Y\cup F}\cdot(\1_n)_{V_i}^{\pm k})$. 

This observation implies that the $n_k$ ($k=1,2,3$) in \eqref{e1} are changeable up to $|{\mu}_{L}(123)|$.  So  
we can suppose that, for all $k=1,2,3$, $n_k$ satisfies 
\begin{equation} \label{e3}
 0\le n_k< |{\mu}_{L}(123)|.  
\end{equation}

Now by \cite{krush} we have, for all $\{i,j,k\}=\{1,2,3\}$, 
 $$ {\mu}_{L}(ijkk)\equiv {\mu}_{\mathrm{cl}(L_0)}(ijkk) + 
{\mu}_{\mathrm{cl}(L_1\cdot L_2\cdot L_3)}(ijkk)\textrm{ mod ${\mu}_{L}(123)$}. $$ 
By Lemma \ref{add}, we have  
$${\mu}_{\mathrm{cl}(L_0)}(ijkk)\equiv 0\textrm{ mod ${\mu}_{L}(123)$}$$
and   
 $$ {\mu}_{\mathrm{cl}(L_1\cdot L_2\cdot L_3)}(ijkk)\equiv\sum_{1\le p\le 3} n_p{\mu}_{\mathrm{cl}(V_p)}(ijkk)
\textrm{ mod ${\mu}_{L}(123)$}. $$ 
  As seen in Remark \ref{remlink}, we have 
${\mu}_{\mathrm{cl}(V_p)}(ijkk)=\delta_{p,k}$. 
It follows that 
\begin{equation} \label{e2}
 {\mu}_{L}(ijkk)\equiv n_k\textrm{ mod ${\mu}_{L}(123)$}.    
\end{equation}

Consider $3$-component Brunnian links $L$ and $L'$ such that 
$\ov\mu_L(123)=\ov\mu_{L'}(123)$ and $\ov\mu_L(ijkk)=\ov\mu_{L'}(ijkk)$ for $(i,j,k)=(1,2,3)$, 
$(1,3,2)$ and $(2,3,1)$.  It follows from (\ref{e1}), (\ref{e2}) and (\ref{e3}) that $L\sim_{C_4} L'$.  
This completes the proof.  
\end{proof}

\subsection{Minimal string link} 

Let $L$ be an $n$-component Brunnian link in $S^3$. 
Denote by $\mathcal{L}(L)$ the set of all $n$-component string links $l$ such that $\mathrm{cl}(l)=L$.  

By Proposition \ref{slcn+1}, for each $l\in \mathcal{L}(L)$ there exists $l'\in SL(n)$ such that $l$ 
is $C_{n+1}$-equivalent to a string link of the form $\prod_{\s \in S_{n-2}} (L_\s)^{\mu_\s(l)}\cdot l'$.  

Put any total order on the set $\mathcal{B}:=\bigcup_{1\le k\le n} \mathcal{B}(k)$ and fix it.  
We denote by $\tau_i$, $i=1,...,m$, the elements of $\mathcal{B}$ according to this total order.  
For all $l\in \mathcal{L}(L)$, $\tau\in \mathcal{B}$, set $\alpha_{\tau}(l):=\mu_{\tau}(l')$.  
For each element $l\in \mathcal{L}(L)$, we can thus define a vector 
  $$ v_l:=( |\alpha_{\tau_1}(l)|,... , |\alpha_{\tau_k}(l)|, ... , |\alpha_{\tau_{m}}(l)|, 
         -\alpha_{\tau_1}(l),... ,-\alpha_{\tau_k}(l), ... , -\alpha_{\tau_{m}}(l)). $$
Set $\mathcal{V}_L=\{ v_l ~|~ l\in \mathcal{L}(L) \}$.  
Define $L_{\mathrm{min}}$ to be the element $l\in \mathcal{L}(L)$ such that $v_l=\min\mathcal{V}_L$ 
(for the natural lexicographical order on $\mathcal{V}$).  
It follows from Proposition \ref{clos} that $L$ is $C_{n+1}$-equivalent to the closure of $L_{\mathrm{min}}$. 
So we have the following.  

\begin{prop}
Two $n$-component Brunnian links $L$ and $L'$ are $C_{n+1}$-equivalent 
if and only if 
$\ov{\mu}_{\sigma}(L)=\ov{\mu}_{\sigma}(L')$ for all $\sigma\in S_{n-1}$ and 
$\min\mathcal{V}_L=\min\mathcal{V}_{L'}$.
\end{prop}

In subsection \ref{S3}, if we take $-|{\mu}_L(123)|/2<n_k<(|{\mu}_L(123)|-1)/2$ 
instead of inequality (\ref{e3}), then we have explicitly $L_{\mathrm{min}}$ for a $3$-component Brunnian link $L$.  
In general, it is a problem to determine $L_{\mathrm{min}}$ from $L$.  

\section{$C_4$-equivalence for links} \label{c4eq}

In this section we prove Theorem \ref{C4} and Proposition \ref{C42}.  
The first subsection provides a lemma which is the main new ingredient for the proofs of these two results.  

\subsection{The index lemma} \label{index}

Let $T$ be a simple $C_k$-tree for an $n$-component link $L$.  
The \emph{index} of $T$ is the collection of all integers $i$ such that $T$ intersects the $i^{th}$ component 
of $L$, counted with multiplicities. For example, a simple $C_3$-tree of index $\{2,3^{(2)},5\}$ for $L$ intersects 
twice component $3$ and once components $2$ and $5$ (and is disjoint from all other components of $L$).  

\begin{lem} \label{ij}
Let $T$ be a simple $C_k$-tree ($k\ge 3$) of index $\{i,j^{(k)}\}$ for an 
$n$-component link $L$, $1\le i\ne j\le n$. Then $L_T\sim_{C_{k+1}} L$.  
\end{lem}

In order to prove this lemma, we need the notion of graph clasper introduced in \cite[\S 8.2]{H}.  
A \emph{graph clasper} is defined as an embedded connected surface 
which is decomposed into leaves, nodes and bands as in Definition \ref{defclasp}, 
but which is not necessarily a disk.  
A graph clasper may contain loops.  
The degree of a graph clasper $G$ is defined as half of the number of nodes and leaves (which coincides with the usual 
degree if $G$ is a tree clasper).  We call a degree $k$ graph clasper a \emph{$C_k$-graph}.  
A $C_k$-graph for a link $L$ is \emph{simple} if each of its leaves intersects $L$ at one point.  

Recall from \cite[\S 8.2]{H} that the STU relation holds for graph clasper.  

\begin{lem}\label{stu}
Let $G_S$, $G_T$ and $G_U$ be three $C_k$-graphs for $\1_n$ which differ only in a small ball 
as depicted in Figure \ref{F06}.  Then 
  $(\1_n)_{G_S}\sim_{C_{k+1}} (\1_n)_{G_T}\cdot (\1_n)_{G_U}$. 
 \begin{figure}[!h]
  \input{stu.pstex_t}
  \caption{The STU relation for $C_k$-graphs.  }\label{F06}
 \end{figure}
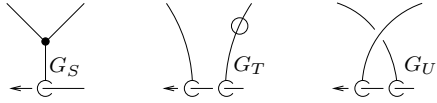
\end{lem}

\noindent It should be noted that, in contrast with the diagram case, this STU 
relation only holds among \emph{connected} claspers.  
Note also that it differs by a sign from the STU relation for unitrivalent diagrams.  

\begin{lem} \label{loop}
Let $C$ be a simple $C_k$-graph for an $n$-component link $L$ in $S^3$, which intersects 
a certain component of $L$ exactly once.  
If $C$ contains a loop (that is, if $C$ is not a $C_k$-tree), then $L_{C}\sim_{C_{k+1}} L$.  
\end{lem}

\begin{proof}
Suppose that $C$ intersect the $i^{th}$ component of $L$ exactly once. 
By \cite{H} and Lemma \ref{cc}, there exists a union $F$ of tree claspers for $\1_n$ and a 
simple $C_k$-tree $G$ for $\1_n$ 
containing a loop and intersecting the $i^{th}$ component once, such that 
 $$ L_C\cong \mathrm{cl}\left( (\1_n)_{F}\cdot (\1_n)_{G} \right). $$
 
Consider the unique leaf $f$ of $G$ intersecting the $i^{th}$ component.  
This leaf $f$ is connected to a loop $\gamma$ of $G$ by a path $P$ of edges and nodes. 
We proceed by induction on the number $n$ of nodes in $P$. 

If $n=0$, that is if $f$ is connected to $\gamma$ by a single edge, apply Lemma \ref{stu} at this edge.  
The result then follows 
from Lemmas \ref{cc} and \ref{twist}, by the arguments similar to those 
in the proof of Proposition \ref{clos}.  

For an arbitrary $n\ge 1$, apply the IHX relation at the edge of $P$ which is incident to $\gamma$.  
By Lemma \ref{asihx},\footnote{Strictly speaking, we cannot apply Lemma \ref{asihx} here, 
as $G$ is not a $C_k$-tree.  However, 
similar relations hold among $C_k$-graphs \cite[\S 8.2]{H}.  } 
we obtain $(\1_n)_{G}\sim_{C_{k+1}} (\1_n)_{G'}\cdot (\1_n)_{G''}$, 
where $G'$ and $G''$ are $C_k$-graphs each of which has a unique leaf intersecting 
the $i^{th}$ component connected to a loop by a path with $(n-1)$ nodes.  
By the induction hypothesis, we thus have 
$(\1_n)_{G'}\sim_{C_{k+1}} \1_n\sim_{C_{k+1}} (\1_n)_{G''}$.  
\end{proof}

\begin{proof}[Proof of Lemma \ref{ij}]
Let $T$ be a simple $C_k$-tree of index $\{i,j^{(k)}\}$ for an 
$n$-component link $L$, $1\le i\ne j\le n$.  
By several applications of Lemmas \ref{stu}, \ref{loop}, \ref{cc} and \ref{twist}, 
one can easily verify that 
$L_T\sim_{C_{k+1}} L_{T'}$, where $T'$ is a simple $C_k$-tree 
of index $\{i,j^{(k)}\}$ for $L$ which 
contains two leaves as depicted in Figure \ref{T}.  
By applying the IHX and STU relations, we have $L_{T'}\sim_{C_{k+1}} L_{T''}$, 
where $T''$ is a $C_k$-graph for $L$ as illustrated in Figure \ref{T}.  
$T''$ clearly satisfies the hypothesis of Lemma \ref{loop}.  
We thus have $L_T\sim_{C_{k+1}} L_{T''}\sim_{C_{k+1}} L$.  
\end{proof}
 \begin{figure}[!h]
  \input{T.pstex_t}
  \caption{}\label{T}
 \end{figure}

\subsection{Proof of Theorem \ref{C4}} \label{c4proof}

We can now prove Theorem \ref{C4}.  
We only need to prove the `if' part of the statement.  

\begin{proof}[Proof of Theorem \ref{C4}]
Let $L$ be a $C_3$-trivial $n$-component link.  
Consider an $n$-component string link $l$ such that its closure is $L$ and $l\sim_{C_3} \1_n$.  
By Lemmas \ref{cc}, \ref{twist} and \ref{asihx}, and the same arguments as those 
used in the proof of Proposition \ref{clos}, we have that 
  $$ l\sim_{C_4} l_0\cdot l_1\cdot l_2\cdot l_3\cdot l_4,  $$ 
where 
 \bi
   \item $l_0=\prod_i (\1_n)_{U_i}$, where $U_i$ is union of simple $C_3$-trees of index $\{i^{(4)}\}$ 
         contained in a regular neighborhood of the $i^{th}$ component of $\1_n$ ; $1\le i\le n$.  
   \item $l_1=\prod_{i<j} \left( (\1_n)_{X_{ij}}\right)^{x_{ij}}$, where $X_{ij}$ is the 
   	 simple $C_3$-tree of index $\{i^{(2)},j^{(2)}\}$ represented in Figure \ref{faim}, and where 
   	 $x_{ij}\in \mathbb{Z}$.  
   \item $l_2=\prod_{i<j ; k} \left( (\1_n)_{Y_{ijk}}\right)^{y_{ijk}}$, where $Y_{ijk}$ is the 
   	 simple $C_3$-tree of index $\{i,j,k^{(2)}\}$ represented in Figure \ref{faim}.
   \item $l_3=\prod_{i\ne j<k<l} \left( (\1_n)_{Z_{ijkl}}\right)^{z_{ijkl}}$, where $Z_{ijkl}$ is the 
   	 simple $C_3$-tree of index $\{i,j,k,l\}$ represented in Figure \ref{faim}.  
   \item $l_4$ is obtained from $\1_n$ by surgery along simple $C_3$-trees with index of the form $\{i,j^{(3)}\}$ 
          ; $1\le i\ne j\le n$.  
 \ei
 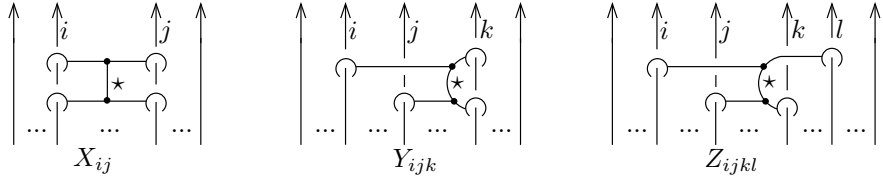
\begin{figure}[!h]
  \input{four.pstex_t}
  \caption{Here $X_{ij}^{-1}$ (resp. $Y_{ijk}^{-1}$, $Z_{ijkl}^{-1}$) is defined as obtained from $X_{ij}$ (resp. 
           $Y_{ijk}$, $Z_{ijkl}$) by a positive half-twist on the edge marked by a $\star$.  }\label{faim}  
 \end{figure}
As an immediate consequence of Lemma \ref{ij}, we thus have 
  $$ L=\mathrm{cl}(l)\sim_{C_4} \mathrm{cl}(l_0\cdot l_1\cdot l_2\cdot l_3). $$ 
It follows from standard computations (see preceding sections) that 
\beQ
\ov{\mu}_L(i,i,j,j) = \mu_{l_1}(i,i,j,j) = x_{ij} & \textrm{for all $1\le i<j\le n$}, \\
\ov{\mu}_L(i,j,k,k) = \mu_{l_2}(i,j,k,k) = y_{ijk} & \textrm{for all $1\le i<j\le n$, $1\le k\le n$}, \\
\ov{\mu}_L(i,j,k,l) = \mu_{l_3}(i,j,k,l) = z_{ijkl} & \textrm{for all $1\le i\ne j<k<l\le n$}.
\eeQ

Now, consider another $C_3$-trivial $n$-component link $L'$, such that $L$ and $L'$ satisfy assertions (1) and (2) in the 
statement of Theorem \ref{C4}.  
By the same construction as above and (1), we have 
 $$ L'\sim_{C_4} \mathrm{cl}(l'_0\cdot l_1\cdot l_2\cdot l_3). $$ 
Here $l'_0=\prod_i (\1_n)_{U'_i}$, where $U'_i$ is union of simple $C_3$-trees of index $\{i^{(4)}\}$ 
contained in a regular neighborhood of the $i^{th}$ component of $\1_n$ ($1\le i\le n$).  
Denote respectively by $(l_0)_i$ and $(l'_0)_i$ the $i^{th}$ component of $l_0$ and $l'_0$.  
By (2) and \cite[Thm. 6.18]{H}, we have $(l_0)_i\sim_{C_4} (l'_0)_i$ for all $i$ in $\{1,...,n\}$.  We thus have 
$l_0\sim_{C_4} l'_0$, which implies the result. 
\end{proof}

\subsection{Proof of Proposition \ref{C42}} \label{c42proof}

It suffices to show that two $2$-component links $L$ and $L'$ which are 
not distinguished by Vassiliev invariants of order $\le 3$
are $C_4$-equivalent (the converse is well-known).  

\begin{proof}[Proof of Proposition \ref{C42}]
By \cite[Thm. 1.5]{M}, $L'$ can be obtained from $L$ by a sequence of surgeries along 
 \be
   \item $C_4$-trees, 
   \item simple $C_3$-trees with index $\{i,j^{(3)}\}$, $\{i,j\}=\{1,2\}$.  
 \ee
By Lemma \ref{ij}, each surgery of type (2) can be achieved by surgery along $C_4$-trees.  
It follows that $L\sim_{C_4} L'$.  
\end{proof}

\end{document}

%% file: asihx.pstex_t
\begin{picture}(0,0)%
\includegraphics{asihx.pstex}%
\end{picture}%
\setlength{\unitlength}{1973sp}%
\begingroup\makeatletter\ifx\SetFigFont\undefined%
\gdef\SetFigFont#1#2#3#4#5{%
  \reset@font\fontsize{#1}{#2pt}%
  \fontfamily{#3}\fontseries{#4}\fontshape{#5}%
  \selectfont}%
\fi\endgroup%
\begin{picture}(8401,1104)(139,-478)
\put(4338,-19){\makebox(0,0)[lb]{\smash{\SetFigFont{8}{9.6}{\rmdefault}{\mddefault}{\updefault}{\color[rgb]{0,0,0}$T_I$}%
}}}
\put(5922,-304){\makebox(0,0)[lb]{\smash{\SetFigFont{8}{9.6}{\rmdefault}{\mddefault}{\updefault}{\color[rgb]{0,0,0}$T_H$}%
}}}
\put(7749,-294){\makebox(0,0)[lb]{\smash{\SetFigFont{8}{9.6}{\rmdefault}{\mddefault}{\updefault}{\color[rgb]{0,0,0}$T_X$}%
}}}
\put(687,-254){\makebox(0,0)[lb]{\smash{\SetFigFont{8}{9.6}{\rmdefault}{\mddefault}{\updefault}{\color[rgb]{0,0,0}$T$}%
}}}
\put(2261,-254){\makebox(0,0)[lb]{\smash{\SetFigFont{8}{9.6}{\rmdefault}{\mddefault}{\updefault}{\color[rgb]{0,0,0}$T'$}%
}}}
\end{picture}

%% file: split.pstex_t
\begin{picture}(0,0)%
\includegraphics{split.pstex}%
\end{picture}%
\setlength{\unitlength}{3947sp}%
\begingroup\makeatletter\ifx\SetFigFont\undefined%
\gdef\SetFigFont#1#2#3#4#5{%
  \reset@font\fontsize{#1}{#2pt}%
  \fontfamily{#3}\fontseries{#4}\fontshape{#5}%
  \selectfont}%
\fi\endgroup%
\begin{picture}(2398,664)(86,-1210)
\put(2386,-1156){\makebox(0,0)[lb]{\smash{\SetFigFont{10}{12.0}{\rmdefault}{\mddefault}{\updefault}{\color[rgb]{0,0,0}$G_2$}%
}}}
\put( 86,-961){\makebox(0,0)[lb]{\smash{\SetFigFont{10}{12.0}{\rmdefault}{\mddefault}{\updefault}{\color[rgb]{0,0,0}$f$}%
}}}
\put(436,-831){\makebox(0,0)[lb]{\smash{\SetFigFont{10}{12.0}{\rmdefault}{\mddefault}{\updefault}{\color[rgb]{0,0,0}$\alpha$}%
}}}
\put(1621,-811){\makebox(0,0)[lb]{\smash{\SetFigFont{10}{12.0}{\rmdefault}{\mddefault}{\updefault}{\color[rgb]{0,0,0}$f_1$}%
}}}
\put(2131,-811){\makebox(0,0)[lb]{\smash{\SetFigFont{10}{12.0}{\rmdefault}{\mddefault}{\updefault}{\color[rgb]{0,0,0}$f_2$}%
}}}
\put(431,-1166){\makebox(0,0)[lb]{\smash{\SetFigFont{10}{12.0}{\rmdefault}{\mddefault}{\updefault}{\color[rgb]{0,0,0}$G$}%
}}}
\put(1876,-1161){\makebox(0,0)[lb]{\smash{\SetFigFont{10}{12.0}{\rmdefault}{\mddefault}{\updefault}{\color[rgb]{0,0,0}$G_1$}%
}}}
\end{picture}

%% file: Tsigm.pstex_t
\begin{picture}(0,0)%
\includegraphics{Tsigm.pstex}%
\end{picture}%
\setlength{\unitlength}{3947sp}%
\begingroup\makeatletter\ifx\SetFigFont\undefined%
\gdef\SetFigFont#1#2#3#4#5{%
  \reset@font\fontsize{#1}{#2pt}%
  \fontfamily{#3}\fontseries{#4}\fontshape{#5}%
  \selectfont}%
\fi\endgroup%
\begin{picture}(2486,1284)(566,-731)
\put(1652,-38){\makebox(0,0)[lb]{\smash{\SetFigFont{10}{12.0}{\rmdefault}{\mddefault}{\updefault}{\color[rgb]{0,0,0}$\sigma$}%
}}}
\put(566,-691){\makebox(0,0)[lb]{\smash{\SetFigFont{10}{12.0}{\rmdefault}{\mddefault}{\updefault}{\color[rgb]{0,0,0}$1$}%
}}}
\put(785,-691){\makebox(0,0)[lb]{\smash{\SetFigFont{10}{12.0}{\rmdefault}{\mddefault}{\updefault}{\color[rgb]{0,0,0}$2$}%
}}}
\put(1846,-191){\makebox(0,0)[lb]{\smash{\SetFigFont{8}{9.6}{\rmdefault}{\mddefault}{\updefault}{\color[rgb]{0,0,0}$n-3$}%
}}}
\put(1841,-466){\makebox(0,0)[lb]{\smash{\SetFigFont{8}{9.6}{\rmdefault}{\mddefault}{\updefault}{\color[rgb]{0,0,0}$n-2$}%
}}}
\put(1004,-691){\makebox(0,0)[lb]{\smash{\SetFigFont{10}{12.0}{\rmdefault}{\mddefault}{\updefault}{\color[rgb]{0,0,0}$n-3$}%
}}}
\put(1414,-691){\makebox(0,0)[lb]{\smash{\SetFigFont{10}{12.0}{\rmdefault}{\mddefault}{\updefault}{\color[rgb]{0,0,0}$n-2$}%
}}}
\put(2213, 64){\rotatebox{270.0}{\makebox(0,0)[lb]{\smash{\SetFigFont{10}{12.0}{\rmdefault}{\mddefault}{\updefault}{\color[rgb]{0,0,0}...}%
}}}}
\put(2643,263){\makebox(0,0)[lb]{\smash{\SetFigFont{10}{12.0}{\rmdefault}{\mddefault}{\updefault}{\color[rgb]{0,0,0}$e$}%
}}}
\put(2891,-686){\makebox(0,0)[lb]{\smash{\SetFigFont{10}{12.0}{\rmdefault}{\mddefault}{\updefault}{\color[rgb]{0,0,0}$n$}%
}}}
\put(2396,-686){\makebox(0,0)[lb]{\smash{\SetFigFont{10}{12.0}{\rmdefault}{\mddefault}{\updefault}{\color[rgb]{0,0,0}$n-1$}%
}}}
\put(1881,379){\makebox(0,0)[lb]{\smash{\SetFigFont{8}{9.6}{\rmdefault}{\mddefault}{\updefault}{\color[rgb]{0,0,0}$1$}%
}}}
\put(1881,104){\makebox(0,0)[lb]{\smash{\SetFigFont{8}{9.6}{\rmdefault}{\mddefault}{\updefault}{\color[rgb]{0,0,0}$2$}%
}}}
\end{picture}

%% file: Tta.pstex_t
\begin{picture}(0,0)%
\includegraphics{Tta.pstex}%
\end{picture}%
\setlength{\unitlength}{3947sp}%
\begingroup\makeatletter\ifx\SetFigFont\undefined%
\gdef\SetFigFont#1#2#3#4#5{%
  \reset@font\fontsize{#1}{#2pt}%
  \fontfamily{#3}\fontseries{#4}\fontshape{#5}%
  \selectfont}%
\fi\endgroup%
\begin{picture}(3851,1548)(181,-996)
\put(1836,-236){\makebox(0,0)[lb]{\smash{\SetFigFont{8}{9.6}{\rmdefault}{\mddefault}{\updefault}{\color[rgb]{0,0,0}$n-2$}%
}}}
\put(1836,-451){\makebox(0,0)[lb]{\smash{\SetFigFont{8}{9.6}{\rmdefault}{\mddefault}{\updefault}{\color[rgb]{0,0,0}$n-1$}%
}}}
\put(1851,379){\makebox(0,0)[lb]{\smash{\SetFigFont{8}{9.6}{\rmdefault}{\mddefault}{\updefault}{\color[rgb]{0,0,0}$1$}%
}}}
\put(1856,149){\makebox(0,0)[lb]{\smash{\SetFigFont{8}{9.6}{\rmdefault}{\mddefault}{\updefault}{\color[rgb]{0,0,0}$2$}%
}}}
\put(925,-952){\makebox(0,0)[lb]{\smash{\SetFigFont{10}{12.0}{\rmdefault}{\mddefault}{\updefault}{\color[rgb]{0,0,0}$k-1$}%
}}}
\put(1626,-38){\makebox(0,0)[lb]{\smash{\SetFigFont{10}{12.0}{\rmdefault}{\mddefault}{\updefault}{\color[rgb]{0,0,0}$\tau$}%
}}}
\put(1506,-95){\rotatebox{270.0}{\makebox(0,0)[lb]{\smash{\SetFigFont{10}{12.0}{\rmdefault}{\mddefault}{\updefault}{\color[rgb]{0,0,0}...}%
}}}}
\put(1506,211){\rotatebox{270.0}{\makebox(0,0)[lb]{\smash{\SetFigFont{10}{12.0}{\rmdefault}{\mddefault}{\updefault}{\color[rgb]{0,0,0}...}%
}}}}
\put(537,-951){\makebox(0,0)[lb]{\smash{\SetFigFont{10}{12.0}{\rmdefault}{\mddefault}{\updefault}{\color[rgb]{0,0,0}$2$}%
}}}
\put(181,-951){\makebox(0,0)[lb]{\smash{\SetFigFont{10}{12.0}{\rmdefault}{\mddefault}{\updefault}{\color[rgb]{0,0,0}$1$}%
}}}
\put(2323,212){\rotatebox{270.0}{\makebox(0,0)[lb]{\smash{\SetFigFont{10}{12.0}{\rmdefault}{\mddefault}{\updefault}{\color[rgb]{0,0,0}...}%
}}}}
\put(2323,-95){\rotatebox{270.0}{\makebox(0,0)[lb]{\smash{\SetFigFont{10}{12.0}{\rmdefault}{\mddefault}{\updefault}{\color[rgb]{0,0,0}...}%
}}}}
\put(2545,434){\makebox(0,0)[lb]{\smash{\SetFigFont{10}{12.0}{\rmdefault}{\mddefault}{\updefault}{\color[rgb]{0,0,0}$e$}%
}}}
\put(3972,-956){\makebox(0,0)[lb]{\smash{\SetFigFont{10}{12.0}{\rmdefault}{\mddefault}{\updefault}{\color[rgb]{0,0,0}$n$}%
}}}
\put(3510,-956){\makebox(0,0)[lb]{\smash{\SetFigFont{10}{12.0}{\rmdefault}{\mddefault}{\updefault}{\color[rgb]{0,0,0}$n-1$}%
}}}
\put(2658,-952){\makebox(0,0)[lb]{\smash{\SetFigFont{10}{12.0}{\rmdefault}{\mddefault}{\updefault}{\color[rgb]{0,0,0}$k$}%
}}}
\put(3045,-952){\makebox(0,0)[lb]{\smash{\SetFigFont{10}{12.0}{\rmdefault}{\mddefault}{\updefault}{\color[rgb]{0,0,0}$k+1$}%
}}}
\end{picture}

%% file: F02.pstex_t
\begin{picture}(0,0)%
\includegraphics{F02.pstex}%
\end{picture}%
\setlength{\unitlength}{3947sp}%
\begingroup\makeatletter\ifx\SetFigFont\undefined%
\gdef\SetFigFont#1#2#3#4#5{%
  \reset@font\fontsize{#1}{#2pt}%
  \fontfamily{#3}\fontseries{#4}\fontshape{#5}%
  \selectfont}%
\fi\endgroup%
\begin{picture}(5123,1538)(128,-1436)
\put(2257,-1396){\makebox(0,0)[lb]{\smash{\SetFigFont{10}{12.0}{\rmdefault}{\mddefault}{\updefault}{\color[rgb]{0,0,0}$1$}%
}}}
\put(2489,-1396){\makebox(0,0)[lb]{\smash{\SetFigFont{10}{12.0}{\rmdefault}{\mddefault}{\updefault}{\color[rgb]{0,0,0}$2$}%
}}}
\put(2718,-1396){\makebox(0,0)[lb]{\smash{\SetFigFont{10}{12.0}{\rmdefault}{\mddefault}{\updefault}{\color[rgb]{0,0,0}$3$}%
}}}
\put(3097,-1396){\makebox(0,0)[lb]{\smash{\SetFigFont{10}{12.0}{\rmdefault}{\mddefault}{\updefault}{\color[rgb]{0,0,0}$4$}%
}}}
\put(2035,-1396){\makebox(0,0)[lb]{\smash{\SetFigFont{10}{12.0}{\rmdefault}{\mddefault}{\updefault}{\color[rgb]{0,0,0}$5$}%
}}}
\put(210,-1396){\makebox(0,0)[lb]{\smash{\SetFigFont{10}{12.0}{\rmdefault}{\mddefault}{\updefault}{\color[rgb]{0,0,0}$1$}%
}}}
\put(442,-1396){\makebox(0,0)[lb]{\smash{\SetFigFont{10}{12.0}{\rmdefault}{\mddefault}{\updefault}{\color[rgb]{0,0,0}$2$}%
}}}
\put(671,-1396){\makebox(0,0)[lb]{\smash{\SetFigFont{10}{12.0}{\rmdefault}{\mddefault}{\updefault}{\color[rgb]{0,0,0}$3$}%
}}}
\put(1061,-1396){\makebox(0,0)[lb]{\smash{\SetFigFont{10}{12.0}{\rmdefault}{\mddefault}{\updefault}{\color[rgb]{0,0,0}$4$}%
}}}
\put(1876,-736){\makebox(0,0)[lb]{\smash{\SetFigFont{10}{12.0}{\rmdefault}{\mddefault}{\updefault}{\color[rgb]{0,0,0}$f$}%
}}}
\put(4126,-436){\makebox(0,0)[lb]{\smash{\SetFigFont{10}{12.0}{\rmdefault}{\mddefault}{\updefault}{\color[rgb]{0,0,0}$f_1$}%
}}}
\put(3901,-1141){\makebox(0,0)[lb]{\smash{\SetFigFont{10}{12.0}{\rmdefault}{\mddefault}{\updefault}{\color[rgb]{0,0,0}$f_2$}%
}}}
\put(4341,-1396){\makebox(0,0)[lb]{\smash{\SetFigFont{10}{12.0}{\rmdefault}{\mddefault}{\updefault}{\color[rgb]{0,0,0}$1$}%
}}}
\put(4573,-1396){\makebox(0,0)[lb]{\smash{\SetFigFont{10}{12.0}{\rmdefault}{\mddefault}{\updefault}{\color[rgb]{0,0,0}$2$}%
}}}
\put(4802,-1396){\makebox(0,0)[lb]{\smash{\SetFigFont{10}{12.0}{\rmdefault}{\mddefault}{\updefault}{\color[rgb]{0,0,0}$3$}%
}}}
\put(4119,-1396){\makebox(0,0)[lb]{\smash{\SetFigFont{10}{12.0}{\rmdefault}{\mddefault}{\updefault}{\color[rgb]{0,0,0}$5$}%
}}}
\put(5251,-1396){\makebox(0,0)[lb]{\smash{\SetFigFont{10}{12.0}{\rmdefault}{\mddefault}{\updefault}{\color[rgb]{0,0,0}$4$}%
}}}
\put(1081,-206){\makebox(0,0)[lb]{\smash{\SetFigFont{10}{12.0}{\rmdefault}{\mddefault}{\updefault}{\color[rgb]{0,0,0}$G_{\eta}$}%
}}}
\put(5216,-191){\makebox(0,0)[lb]{\smash{\SetFigFont{10}{12.0}{\rmdefault}{\mddefault}{\updefault}{\color[rgb]{0,0,0}$G^1_{\eta}$}%
}}}
\put(5216,-916){\makebox(0,0)[lb]{\smash{\SetFigFont{10}{12.0}{\rmdefault}{\mddefault}{\updefault}{\color[rgb]{0,0,0}$G^2_{\eta}$}%
}}}
\put(3116,-206){\makebox(0,0)[lb]{\smash{\SetFigFont{10}{12.0}{\rmdefault}{\mddefault}{\updefault}{\color[rgb]{0,0,0}$\tilde{G}_{\eta}$}%
}}}
\end{picture}

%% file: F01.pstex_t
\begin{picture}(0,0)%
\includegraphics{F01.pstex}%
\end{picture}%
\setlength{\unitlength}{3947sp}%
\begingroup\makeatletter\ifx\SetFigFont\undefined%
\gdef\SetFigFont#1#2#3#4#5{%
  \reset@font\fontsize{#1}{#2pt}%
  \fontfamily{#3}\fontseries{#4}\fontshape{#5}%
  \selectfont}%
\fi\endgroup%
\begin{picture}(4860,1912)(86,-1511)
\put(1005,-1471){\makebox(0,0)[lb]{\smash{\SetFigFont{10}{12.0}{\rmdefault}{\mddefault}{\updefault}{\color[rgb]{0,0,0}$4$}%
}}}
\put(168,-1471){\makebox(0,0)[lb]{\smash{\SetFigFont{10}{12.0}{\rmdefault}{\mddefault}{\updefault}{\color[rgb]{0,0,0}$1$}%
}}}
\put(400,-1471){\makebox(0,0)[lb]{\smash{\SetFigFont{10}{12.0}{\rmdefault}{\mddefault}{\updefault}{\color[rgb]{0,0,0}$2$}%
}}}
\put(629,-1471){\makebox(0,0)[lb]{\smash{\SetFigFont{10}{12.0}{\rmdefault}{\mddefault}{\updefault}{\color[rgb]{0,0,0}$3$}%
}}}
\put(2709,-1471){\makebox(0,0)[lb]{\smash{\SetFigFont{10}{12.0}{\rmdefault}{\mddefault}{\updefault}{\color[rgb]{0,0,0}$4$}%
}}}
\put(2929,-1471){\makebox(0,0)[lb]{\smash{\SetFigFont{10}{12.0}{\rmdefault}{\mddefault}{\updefault}{\color[rgb]{0,0,0}$5$}%
}}}
\put(1869,-1471){\makebox(0,0)[lb]{\smash{\SetFigFont{10}{12.0}{\rmdefault}{\mddefault}{\updefault}{\color[rgb]{0,0,0}$1$}%
}}}
\put(2101,-1471){\makebox(0,0)[lb]{\smash{\SetFigFont{10}{12.0}{\rmdefault}{\mddefault}{\updefault}{\color[rgb]{0,0,0}$2$}%
}}}
\put(2330,-1471){\makebox(0,0)[lb]{\smash{\SetFigFont{10}{12.0}{\rmdefault}{\mddefault}{\updefault}{\color[rgb]{0,0,0}$3$}%
}}}
\put(3763,-1471){\makebox(0,0)[lb]{\smash{\SetFigFont{10}{12.0}{\rmdefault}{\mddefault}{\updefault}{\color[rgb]{0,0,0}$1$}%
}}}
\put(3995,-1471){\makebox(0,0)[lb]{\smash{\SetFigFont{10}{12.0}{\rmdefault}{\mddefault}{\updefault}{\color[rgb]{0,0,0}$2$}%
}}}
\put(4224,-1471){\makebox(0,0)[lb]{\smash{\SetFigFont{10}{12.0}{\rmdefault}{\mddefault}{\updefault}{\color[rgb]{0,0,0}$3$}%
}}}
\put(4558,-1471){\makebox(0,0)[lb]{\smash{\SetFigFont{10}{12.0}{\rmdefault}{\mddefault}{\updefault}{\color[rgb]{0,0,0}$4$}%
}}}
\put(4902,-1471){\makebox(0,0)[lb]{\smash{\SetFigFont{10}{12.0}{\rmdefault}{\mddefault}{\updefault}{\color[rgb]{0,0,0}$5$}%
}}}
\put(2952,-661){\makebox(0,0)[lb]{\smash{\SetFigFont{10}{12.0}{\rmdefault}{\mddefault}{\updefault}{\color[rgb]{0,0,0}$\tilde{G}_{\eta}$}%
}}}
\put(1024,-706){\makebox(0,0)[lb]{\smash{\SetFigFont{10}{12.0}{\rmdefault}{\mddefault}{\updefault}{\color[rgb]{0,0,0}$G_{\eta}$}%
}}}
\put(4941,138){\makebox(0,0)[lb]{\smash{\SetFigFont{10}{12.0}{\rmdefault}{\mddefault}{\updefault}{\color[rgb]{0,0,0}$G^1_{\eta}$}%
}}}
\put(4938,-297){\makebox(0,0)[lb]{\smash{\SetFigFont{10}{12.0}{\rmdefault}{\mddefault}{\updefault}{\color[rgb]{0,0,0}$G^2_{\eta}$}%
}}}
\put(4940,-781){\makebox(0,0)[lb]{\smash{\SetFigFont{10}{12.0}{\rmdefault}{\mddefault}{\updefault}{\color[rgb]{0,0,0}$G^3_{\eta}$}%
}}}
\put(4933,-1160){\makebox(0,0)[lb]{\smash{\SetFigFont{10}{12.0}{\rmdefault}{\mddefault}{\updefault}{\color[rgb]{0,0,0}$G^4_{\eta}$}%
}}}
\end{picture}

%% file: F03.pstex_t
\begin{picture}(0,0)%
\includegraphics{F03.pstex}%
\end{picture}%
\setlength{\unitlength}{3947sp}%
\begingroup\makeatletter\ifx\SetFigFont\undefined%
\gdef\SetFigFont#1#2#3#4#5{%
  \reset@font\fontsize{#1}{#2pt}%
  \fontfamily{#3}\fontseries{#4}\fontshape{#5}%
  \selectfont}%
\fi\endgroup%
\begin{picture}(5537,1614)(128,-1512)
\put(210,-1471){\makebox(0,0)[lb]{\smash{\SetFigFont{10}{12.0}{\rmdefault}{\mddefault}{\updefault}{\color[rgb]{0,0,0}$1$}%
}}}
\put(442,-1471){\makebox(0,0)[lb]{\smash{\SetFigFont{10}{12.0}{\rmdefault}{\mddefault}{\updefault}{\color[rgb]{0,0,0}$2$}%
}}}
\put(671,-1471){\makebox(0,0)[lb]{\smash{\SetFigFont{10}{12.0}{\rmdefault}{\mddefault}{\updefault}{\color[rgb]{0,0,0}$3$}%
}}}
\put(1061,-1471){\makebox(0,0)[lb]{\smash{\SetFigFont{10}{12.0}{\rmdefault}{\mddefault}{\updefault}{\color[rgb]{0,0,0}$4$}%
}}}
\put(1931,-1472){\makebox(0,0)[lb]{\smash{\SetFigFont{10}{12.0}{\rmdefault}{\mddefault}{\updefault}{\color[rgb]{0,0,0}$1$}%
}}}
\put(2163,-1472){\makebox(0,0)[lb]{\smash{\SetFigFont{10}{12.0}{\rmdefault}{\mddefault}{\updefault}{\color[rgb]{0,0,0}$2$}%
}}}
\put(2392,-1472){\makebox(0,0)[lb]{\smash{\SetFigFont{10}{12.0}{\rmdefault}{\mddefault}{\updefault}{\color[rgb]{0,0,0}$3$}%
}}}
\put(2726,-1472){\makebox(0,0)[lb]{\smash{\SetFigFont{10}{12.0}{\rmdefault}{\mddefault}{\updefault}{\color[rgb]{0,0,0}$4$}%
}}}
\put(3070,-1472){\makebox(0,0)[lb]{\smash{\SetFigFont{10}{12.0}{\rmdefault}{\mddefault}{\updefault}{\color[rgb]{0,0,0}$5$}%
}}}
\put(4148,-1472){\makebox(0,0)[lb]{\smash{\SetFigFont{10}{12.0}{\rmdefault}{\mddefault}{\updefault}{\color[rgb]{0,0,0}$1$}%
}}}
\put(4380,-1472){\makebox(0,0)[lb]{\smash{\SetFigFont{10}{12.0}{\rmdefault}{\mddefault}{\updefault}{\color[rgb]{0,0,0}$2$}%
}}}
\put(4609,-1472){\makebox(0,0)[lb]{\smash{\SetFigFont{10}{12.0}{\rmdefault}{\mddefault}{\updefault}{\color[rgb]{0,0,0}$3$}%
}}}
\put(5381,-1472){\makebox(0,0)[lb]{\smash{\SetFigFont{10}{12.0}{\rmdefault}{\mddefault}{\updefault}{\color[rgb]{0,0,0}$4$}%
}}}
\put(5623,-1472){\makebox(0,0)[lb]{\smash{\SetFigFont{10}{12.0}{\rmdefault}{\mddefault}{\updefault}{\color[rgb]{0,0,0}$5$}%
}}}
\put(3076,-1103){\makebox(0,0)[lb]{\smash{\SetFigFont{10}{12.0}{\rmdefault}{\mddefault}{\updefault}{\color[rgb]{0,0,0}$G^2_{\eta \rho}$}%
}}}
\put(3075,-277){\makebox(0,0)[lb]{\smash{\SetFigFont{10}{12.0}{\rmdefault}{\mddefault}{\updefault}{\color[rgb]{0,0,0}$G^1_{\eta \rho}$}%
}}}
\put(1065,-889){\makebox(0,0)[lb]{\smash{\SetFigFont{10}{12.0}{\rmdefault}{\mddefault}{\updefault}{\color[rgb]{0,0,0}$G_{\eta \rho}$}%
}}}
\put(5630,-412){\makebox(0,0)[lb]{\smash{\SetFigFont{10}{12.0}{\rmdefault}{\mddefault}{\updefault}{\color[rgb]{0,0,0}$T^1_{\eta}$}%
}}}
\put(5646,-1023){\makebox(0,0)[lb]{\smash{\SetFigFont{10}{12.0}{\rmdefault}{\mddefault}{\updefault}{\color[rgb]{0,0,0}$T^2_{\eta}$}%
}}}
\end{picture}

%% file: link.pstex_t
\begin{picture}(0,0)%
\includegraphics{link.pstex}%
\end{picture}%
\setlength{\unitlength}{2763sp}%
\begingroup\makeatletter\ifx\SetFigFont\undefined%
\gdef\SetFigFont#1#2#3#4#5{%
  \reset@font\fontsize{#1}{#2pt}%
  \fontfamily{#3}\fontseries{#4}\fontshape{#5}%
  \selectfont}%
\fi\endgroup%
\begin{picture}(3980,1810)(589,-1476)
\put(602,-1191){\makebox(0,0)[lb]{\smash{\SetFigFont{7}{8.4}{\rmdefault}{\mddefault}{\updefault}{\color[rgb]{0,0,0}$\tau(1)$}%
}}}
\put(1707,-1296){\makebox(0,0)[lb]{\smash{\SetFigFont{7}{8.4}{\rmdefault}{\mddefault}{\updefault}{\color[rgb]{0,0,0}$\tau(2)$}%
}}}
\put(2502,-1296){\makebox(0,0)[lb]{\smash{\SetFigFont{7}{8.4}{\rmdefault}{\mddefault}{\updefault}{\color[rgb]{0,0,0}$\tau(3)$}%
}}}
\put(2607,174){\makebox(0,0)[lb]{\smash{\SetFigFont{7}{8.4}{\rmdefault}{\mddefault}{\updefault}{\color[rgb]{0,0,0}$k$}%
}}}
\put(4567,-1081){\makebox(0,0)[lb]{\smash{\SetFigFont{7}{8.4}{\rmdefault}{\mddefault}{\updefault}{\color[rgb]{0,0,0}$\tau(n-1)$}%
}}}
\end{picture}

%% file: Ktau.pstex_t
\begin{picture}(0,0)%
\includegraphics{Ktau.pstex}%
\end{picture}%
\setlength{\unitlength}{3947sp}%
\begingroup\makeatletter\ifx\SetFigFont\undefined%
\gdef\SetFigFont#1#2#3#4#5{%
  \reset@font\fontsize{#1}{#2pt}%
  \fontfamily{#3}\fontseries{#4}\fontshape{#5}%
  \selectfont}%
\fi\endgroup%
\begin{picture}(3396,1377)(461,-1218)
\put(614,-1185){\makebox(0,0)[lb]{\smash{\SetFigFont{8}{9.6}{\rmdefault}{\mddefault}{\updefault}{\color[rgb]{0,0,0}$1$}%
}}}
\put(874,-1185){\makebox(0,0)[lb]{\smash{\SetFigFont{8}{9.6}{\rmdefault}{\mddefault}{\updefault}{\color[rgb]{0,0,0}$2$}%
}}}
\put(1132,-1185){\makebox(0,0)[lb]{\smash{\SetFigFont{8}{9.6}{\rmdefault}{\mddefault}{\updefault}{\color[rgb]{0,0,0}$3$}%
}}}
\put(1386,-1185){\makebox(0,0)[lb]{\smash{\SetFigFont{8}{9.6}{\rmdefault}{\mddefault}{\updefault}{\color[rgb]{0,0,0}$4$}%
}}}
\put(897,-550){\makebox(0,0)[lb]{\smash{\SetFigFont{8}{9.6}{\rmdefault}{\mddefault}{\updefault}{\color[rgb]{0,0,0}$L$}%
}}}
\put(3258,-550){\makebox(0,0)[lb]{\smash{\SetFigFont{8}{9.6}{\rmdefault}{\mddefault}{\updefault}{\color[rgb]{0,0,0}$L$}%
}}}
\put(3777,-34){\makebox(0,0)[lb]{\smash{\SetFigFont{8}{9.6}{\rmdefault}{\mddefault}{\updefault}{\color[rgb]{0,0,0}$K_{\tau}(L)$}%
}}}
\end{picture}

%% file: F04.pstex_t
\begin{picture}(0,0)%
\includegraphics{F04.pstex}%
\end{picture}%
\setlength{\unitlength}{3947sp}%
\begingroup\makeatletter\ifx\SetFigFont\undefined%
\gdef\SetFigFont#1#2#3#4#5{%
  \reset@font\fontsize{#1}{#2pt}%
  \fontfamily{#3}\fontseries{#4}\fontshape{#5}%
  \selectfont}%
\fi\endgroup%
\begin{picture}(5809,1613)(128,-1511)
\put(210,-1471){\makebox(0,0)[lb]{\smash{\SetFigFont{10}{12.0}{\rmdefault}{\mddefault}{\updefault}{\color[rgb]{0,0,0}$1$}%
}}}
\put(442,-1471){\makebox(0,0)[lb]{\smash{\SetFigFont{10}{12.0}{\rmdefault}{\mddefault}{\updefault}{\color[rgb]{0,0,0}$2$}%
}}}
\put(671,-1471){\makebox(0,0)[lb]{\smash{\SetFigFont{10}{12.0}{\rmdefault}{\mddefault}{\updefault}{\color[rgb]{0,0,0}$3$}%
}}}
\put(4524,-1471){\makebox(0,0)[lb]{\smash{\SetFigFont{10}{12.0}{\rmdefault}{\mddefault}{\updefault}{\color[rgb]{0,0,0}$1$}%
}}}
\put(4756,-1471){\makebox(0,0)[lb]{\smash{\SetFigFont{10}{12.0}{\rmdefault}{\mddefault}{\updefault}{\color[rgb]{0,0,0}$2$}%
}}}
\put(4985,-1471){\makebox(0,0)[lb]{\smash{\SetFigFont{10}{12.0}{\rmdefault}{\mddefault}{\updefault}{\color[rgb]{0,0,0}$3$}%
}}}
\put(1098,-1108){\makebox(0,0)[lb]{\smash{\SetFigFont{10}{12.0}{\rmdefault}{\mddefault}{\updefault}{\color[rgb]{0,0,0}$f$}%
}}}
\put(1060,-754){\makebox(0,0)[lb]{\smash{\SetFigFont{10}{12.0}{\rmdefault}{\mddefault}{\updefault}{\color[rgb]{0,0,0}$4$}%
}}}
\put(5803,-1315){\makebox(0,0)[lb]{\smash{\SetFigFont{10}{12.0}{\rmdefault}{\mddefault}{\updefault}{\color[rgb]{0,0,0}$4$}%
}}}
\put(676,-361){\makebox(0,0)[lb]{\smash{\SetFigFont{10}{12.0}{\rmdefault}{\mddefault}{\updefault}{\color[rgb]{0,0,0}$G$}%
}}}
\put(2275,-1471){\makebox(0,0)[lb]{\smash{\SetFigFont{10}{12.0}{\rmdefault}{\mddefault}{\updefault}{\color[rgb]{0,0,0}$1$}%
}}}
\put(2507,-1471){\makebox(0,0)[lb]{\smash{\SetFigFont{10}{12.0}{\rmdefault}{\mddefault}{\updefault}{\color[rgb]{0,0,0}$2$}%
}}}
\put(2736,-1471){\makebox(0,0)[lb]{\smash{\SetFigFont{10}{12.0}{\rmdefault}{\mddefault}{\updefault}{\color[rgb]{0,0,0}$3$}%
}}}
\put(3315,-1039){\makebox(0,0)[lb]{\smash{\SetFigFont{10}{12.0}{\rmdefault}{\mddefault}{\updefault}{\color[rgb]{0,0,0}$4$}%
}}}
\put(2915,-1111){\makebox(0,0)[lb]{\smash{\SetFigFont{10}{12.0}{\rmdefault}{\mddefault}{\updefault}{\color[rgb]{0,0,0}$G'$}%
}}}
\put(5026,-361){\makebox(0,0)[lb]{\smash{\SetFigFont{10}{12.0}{\rmdefault}{\mddefault}{\updefault}{\color[rgb]{0,0,0}$G''$}%
}}}
\end{picture}

%% file: bor.pstex_t
\begin{picture}(0,0)%
\includegraphics{bor.pstex}%
\end{picture}%
\setlength{\unitlength}{3947sp}%
\begingroup\makeatletter\ifx\SetFigFont\undefined%
\gdef\SetFigFont#1#2#3#4#5{%
  \reset@font\fontsize{#1}{#2pt}%
  \fontfamily{#3}\fontseries{#4}\fontshape{#5}%
  \selectfont}%
\fi\endgroup%
\begin{picture}(5910,931)(83,-240)
\put(2816,-186){\makebox(0,0)[lb]{\smash{\SetFigFont{10}{12.0}{\rmdefault}{\mddefault}{\updefault}{\color[rgb]{0,0,0}$V_1$}%
}}}
\put(4496,-191){\makebox(0,0)[lb]{\smash{\SetFigFont{10}{12.0}{\rmdefault}{\mddefault}{\updefault}{\color[rgb]{0,0,0}$V_2$}%
}}}
\put(676,-181){\makebox(0,0)[lb]{\smash{\SetFigFont{10}{12.0}{\rmdefault}{\mddefault}{\updefault}{\color[rgb]{0,0,0}$B$}%
}}}
\put(5531,-191){\makebox(0,0)[lb]{\smash{\SetFigFont{10}{12.0}{\rmdefault}{\mddefault}{\updefault}{\color[rgb]{0,0,0}$V_3$}%
}}}
\put(5736,134){\makebox(0,0)[lb]{\smash{\SetFigFont{10}{12.0}{\rmdefault}{\mddefault}{\updefault}{\color[rgb]{0,0,0}$\star$}%
}}}
\put(676,449){\makebox(0,0)[lb]{\smash{\SetFigFont{10}{12.0}{\rmdefault}{\mddefault}{\updefault}{\color[rgb]{0,0,0}$\star$}%
}}}
\put(1841,236){\makebox(0,0)[lb]{\smash{\SetFigFont{10}{12.0}{\rmdefault}{\mddefault}{\updefault}{\color[rgb]{0,0,0}$\star$}%
}}}
\put(3835,237){\makebox(0,0)[lb]{\smash{\SetFigFont{10}{12.0}{\rmdefault}{\mddefault}{\updefault}{\color[rgb]{0,0,0}$\star$}%
}}}
\end{picture}

%% file: isotopy.pstex_t
\begin{picture}(0,0)%
\includegraphics{isotopy.pstex}%
\end{picture}%
\setlength{\unitlength}{3947sp}%
\begingroup\makeatletter\ifx\SetFigFont\undefined%
\gdef\SetFigFont#1#2#3#4#5{%
  \reset@font\fontsize{#1}{#2pt}%
  \fontfamily{#3}\fontseries{#4}\fontshape{#5}%
  \selectfont}%
\fi\endgroup%
\begin{picture}(5694,1476)(5,-638)
\put(5050,-586){\makebox(0,0)[lb]{\smash{\SetFigFont{10}{12.0}{\rmdefault}{\mddefault}{\updefault}{\color[rgb]{0,0,0}$j$}%
}}}
\put(5040,429){\makebox(0,0)[lb]{\smash{\SetFigFont{10}{12.0}{\rmdefault}{\mddefault}{\updefault}{\color[rgb]{0,0,0}$i$}%
}}}
\put(4695, 54){\makebox(0,0)[lb]{\smash{\SetFigFont{10}{12.0}{\rmdefault}{\mddefault}{\updefault}{\color[rgb]{0,0,0}$Y$}%
}}}
\put(5648,279){\makebox(0,0)[lb]{\smash{\SetFigFont{10}{12.0}{\rmdefault}{\mddefault}{\updefault}{\color[rgb]{0,0,0}$k$}%
}}}
\put(2589,-583){\makebox(0,0)[lb]{\smash{\SetFigFont{10}{12.0}{\rmdefault}{\mddefault}{\updefault}{\color[rgb]{0,0,0}$j$}%
}}}
\put(2259,664){\makebox(0,0)[lb]{\smash{\SetFigFont{10}{12.0}{\rmdefault}{\mddefault}{\updefault}{\color[rgb]{0,0,0}$i$}%
}}}
\put(2936,390){\makebox(0,0)[lb]{\smash{\SetFigFont{10}{12.0}{\rmdefault}{\mddefault}{\updefault}{\color[rgb]{0,0,0}$Y'$}%
}}}
\put(3354,310){\makebox(0,0)[lb]{\smash{\SetFigFont{10}{12.0}{\rmdefault}{\mddefault}{\updefault}{\color[rgb]{0,0,0}$k$}%
}}}
\put(473,422){\makebox(0,0)[lb]{\smash{\SetFigFont{10}{12.0}{\rmdefault}{\mddefault}{\updefault}{\color[rgb]{0,0,0}$i$}%
}}}
\put(483,-593){\makebox(0,0)[lb]{\smash{\SetFigFont{10}{12.0}{\rmdefault}{\mddefault}{\updefault}{\color[rgb]{0,0,0}$j$}%
}}}
\put(128, 47){\makebox(0,0)[lb]{\smash{\SetFigFont{10}{12.0}{\rmdefault}{\mddefault}{\updefault}{\color[rgb]{0,0,0}$Y$}%
}}}
\put(999,332){\makebox(0,0)[lb]{\smash{\SetFigFont{10}{12.0}{\rmdefault}{\mddefault}{\updefault}{\color[rgb]{0,0,0}$k$}%
}}}
\put(3645,105){\makebox(0,0)[lb]{\smash{\SetFigFont{10}{12.0}{\rmdefault}{\mddefault}{\updefault}{\color[rgb]{0,0,0}$C_3$-moves}%
}}}
\end{picture}

%% file: stu.pstex_t
\begin{picture}(0,0)%
\includegraphics{stu.pstex}%
\end{picture}%
\setlength{\unitlength}{1973sp}%
\begingroup\makeatletter\ifx\SetFigFont\undefined%
\gdef\SetFigFont#1#2#3#4#5{%
  \reset@font\fontsize{#1}{#2pt}%
  \fontfamily{#3}\fontseries{#4}\fontshape{#5}%
  \selectfont}%
\fi\endgroup%
\begin{picture}(5238,1192)(139,-565)
\put(687,-254){\makebox(0,0)[lb]{\smash{\SetFigFont{8}{9.6}{\rmdefault}{\mddefault}{\updefault}{\color[rgb]{0,0,0}$G_S$}%
}}}
\put(5144,-234){\makebox(0,0)[lb]{\smash{\SetFigFont{8}{9.6}{\rmdefault}{\mddefault}{\updefault}{\color[rgb]{0,0,0}$G_U$}%
}}}
\put(3002,-234){\makebox(0,0)[lb]{\smash{\SetFigFont{8}{9.6}{\rmdefault}{\mddefault}{\updefault}{\color[rgb]{0,0,0}$G_T$}%
}}}
\end{picture}

%% file: T.pstex_t
\begin{picture}(0,0)%
\includegraphics{T.pstex}%
\end{picture}%
\setlength{\unitlength}{1973sp}%
\begingroup\makeatletter\ifx\SetFigFont\undefined%
\gdef\SetFigFont#1#2#3#4#5{%
  \reset@font\fontsize{#1}{#2pt}%
  \fontfamily{#3}\fontseries{#4}\fontshape{#5}%
  \selectfont}%
\fi\endgroup%
\begin{picture}(4890,1689)(301,-1824)
\put(1027,-862){\makebox(0,0)[lb]{\smash{\SetFigFont{6}{7.2}{\rmdefault}{\mddefault}{\updefault}{\color[rgb]{0,0,0}$T'$}%
}}}
\put(4295,-862){\makebox(0,0)[lb]{\smash{\SetFigFont{6}{7.2}{\rmdefault}{\mddefault}{\updefault}{\color[rgb]{0,0,0}$T''$}%
}}}
\put(301,-1770){\makebox(0,0)[lb]{\smash{\SetFigFont{6}{7.2}{\rmdefault}{\mddefault}{\updefault}{\color[rgb]{0,0,0}$j$}%
}}}
\put(3884,-1770){\makebox(0,0)[lb]{\smash{\SetFigFont{6}{7.2}{\rmdefault}{\mddefault}{\updefault}{\color[rgb]{0,0,0}$j$}%
}}}
\end{picture}

%% file: four.pstex_t
\begin{picture}(0,0)%
\includegraphics{four.pstex}%
\end{picture}%
\setlength{\unitlength}{3947sp}%
\begingroup\makeatletter\ifx\SetFigFont\undefined%
\gdef\SetFigFont#1#2#3#4#5{%
  \reset@font\fontsize{#1}{#2pt}%
  \fontfamily{#3}\fontseries{#4}\fontshape{#5}%
  \selectfont}%
\fi\endgroup%
\begin{picture}(5494,1101)(131,-410)
\put(5394,-141){\makebox(0,0)[lb]{\smash{\SetFigFont{10}{12.0}{\rmdefault}{\mddefault}{\updefault}{\color[rgb]{0,0,0}...}%
}}}
\put(5123,-141){\makebox(0,0)[lb]{\smash{\SetFigFont{10}{12.0}{\rmdefault}{\mddefault}{\updefault}{\color[rgb]{0,0,0}...}%
}}}
\put(4019,-141){\makebox(0,0)[lb]{\smash{\SetFigFont{10}{12.0}{\rmdefault}{\mddefault}{\updefault}{\color[rgb]{0,0,0}...}%
}}}
\put(4231,444){\makebox(0,0)[lb]{\smash{\SetFigFont{10}{12.0}{\rmdefault}{\mddefault}{\updefault}{\color[rgb]{0,0,0}$i$}%
}}}
\put(4349,-141){\makebox(0,0)[lb]{\smash{\SetFigFont{10}{12.0}{\rmdefault}{\mddefault}{\updefault}{\color[rgb]{0,0,0}...}%
}}}
\put(4722,-141){\makebox(0,0)[lb]{\smash{\SetFigFont{10}{12.0}{\rmdefault}{\mddefault}{\updefault}{\color[rgb]{0,0,0}...}%
}}}
\put(4605,447){\makebox(0,0)[lb]{\smash{\SetFigFont{10}{12.0}{\rmdefault}{\mddefault}{\updefault}{\color[rgb]{0,0,0}$j$}%
}}}
\put(5055,447){\makebox(0,0)[lb]{\smash{\SetFigFont{10}{12.0}{\rmdefault}{\mddefault}{\updefault}{\color[rgb]{0,0,0}$k$}%
}}}
\put(5337,447){\makebox(0,0)[lb]{\smash{\SetFigFont{10}{12.0}{\rmdefault}{\mddefault}{\updefault}{\color[rgb]{0,0,0}$l$}%
}}}
\put(464,444){\makebox(0,0)[lb]{\smash{\SetFigFont{10}{12.0}{\rmdefault}{\mddefault}{\updefault}{\color[rgb]{0,0,0}$i$}%
}}}
\put(1091,447){\makebox(0,0)[lb]{\smash{\SetFigFont{10}{12.0}{\rmdefault}{\mddefault}{\updefault}{\color[rgb]{0,0,0}$j$}%
}}}
\put(707,-141){\makebox(0,0)[lb]{\smash{\SetFigFont{10}{12.0}{\rmdefault}{\mddefault}{\updefault}{\color[rgb]{0,0,0}...}%
}}}
\put(250,-141){\makebox(0,0)[lb]{\smash{\SetFigFont{10}{12.0}{\rmdefault}{\mddefault}{\updefault}{\color[rgb]{0,0,0}...}%
}}}
\put(1159,-141){\makebox(0,0)[lb]{\smash{\SetFigFont{10}{12.0}{\rmdefault}{\mddefault}{\updefault}{\color[rgb]{0,0,0}...}%
}}}
\put(3167,-141){\makebox(0,0)[lb]{\smash{\SetFigFont{10}{12.0}{\rmdefault}{\mddefault}{\updefault}{\color[rgb]{0,0,0}...}%
}}}
\put(2063,-141){\makebox(0,0)[lb]{\smash{\SetFigFont{10}{12.0}{\rmdefault}{\mddefault}{\updefault}{\color[rgb]{0,0,0}...}%
}}}
\put(2275,444){\makebox(0,0)[lb]{\smash{\SetFigFont{10}{12.0}{\rmdefault}{\mddefault}{\updefault}{\color[rgb]{0,0,0}$i$}%
}}}
\put(2393,-141){\makebox(0,0)[lb]{\smash{\SetFigFont{10}{12.0}{\rmdefault}{\mddefault}{\updefault}{\color[rgb]{0,0,0}...}%
}}}
\put(2766,-141){\makebox(0,0)[lb]{\smash{\SetFigFont{10}{12.0}{\rmdefault}{\mddefault}{\updefault}{\color[rgb]{0,0,0}...}%
}}}
\put(2649,447){\makebox(0,0)[lb]{\smash{\SetFigFont{10}{12.0}{\rmdefault}{\mddefault}{\updefault}{\color[rgb]{0,0,0}$j$}%
}}}
\put(3099,447){\makebox(0,0)[lb]{\smash{\SetFigFont{10}{12.0}{\rmdefault}{\mddefault}{\updefault}{\color[rgb]{0,0,0}$k$}%
}}}
\put(4512,-361){\makebox(0,0)[lb]{\smash{\SetFigFont{10}{12.0}{\rmdefault}{\mddefault}{\updefault}{\color[rgb]{0,0,0}$Z_{ijkl}$}%
}}}
\put(2556,-361){\makebox(0,0)[lb]{\smash{\SetFigFont{10}{12.0}{\rmdefault}{\mddefault}{\updefault}{\color[rgb]{0,0,0}$Y_{ijk}$}%
}}}
\put(548,-357){\makebox(0,0)[lb]{\smash{\SetFigFont{10}{12.0}{\rmdefault}{\mddefault}{\updefault}{\color[rgb]{0,0,0}$X_{ij}$}%
}}}
\put(788,126){\makebox(0,0)[lb]{\smash{\SetFigFont{10}{12.0}{\rmdefault}{\mddefault}{\updefault}{\color[rgb]{0,0,0}$\star$}%
}}}
\put(4876,137){\makebox(0,0)[lb]{\smash{\SetFigFont{10}{12.0}{\rmdefault}{\mddefault}{\updefault}{\color[rgb]{0,0,0}$\star$}%
}}}
\put(2919,127){\makebox(0,0)[lb]{\smash{\SetFigFont{10}{12.0}{\rmdefault}{\mddefault}{\updefault}{\color[rgb]{0,0,0}$\star$}%
}}}
\end{picture}